\documentclass{amsart}
\usepackage{amssymb,amsmath}
\usepackage[mathscr]{euscript}
\input{epsf}

\usepackage[cp1251]{inputenc}
\usepackage[greek,  english]{babel}

\usepackage{teubner}
\usepackage{bbding}

\newcounter{sec}

\newcounter{punct}[sec]

\def\punct{\refstepcounter{punct}{\arabic{sec}.\arabic{punct}.  }}

\newtheorem{theorem}{Theorem}[sec]
\newtheorem{proposition}[theorem]{Proposition}

\newtheorem{lemma}[theorem]{Lemma}

\def\COUNTERS{\addtocounter{sec}{1}
              \setcounter{punct}{0}
          \setcounter{equation}{0}
          \setcounter{theorem}{0}
          }
          
          \def\sm{\smallskip}

\begin{document}

\newcommand{\supp}{\mathop {\mathrm {supp}}\nolimits}
\newcommand{\rk}{\mathop {\mathrm {rk}}\nolimits}
\newcommand{\Aut}{\mathop {\mathrm {Aut}}\nolimits}
\newcommand{\Ob}{\mathop {\mathrm {Ob}}\nolimits}
\newcommand{\Out}{\mathop {\mathrm {Out}}\nolimits}
\renewcommand{\Re}{\mathop {\mathrm {Re}}\nolimits}
\newcommand{\Inn}{\mathop {\mathrm {Inn}}\nolimits}
\newcommand{\Char}{\mathop {\mathrm {Char}}\nolimits}
\newcommand{\ch}{\cosh}
\newcommand{\sh}{\sinh}
\newcommand{\Sp}{\mathop {\mathrm {Sp}}\nolimits}
\newcommand{\SOS}{\mathop {\mathrm {SO^*}}\nolimits}
\newcommand{\Ams}{\mathop {\mathrm {Ams}}\nolimits}
\newcommand{\Gms}{\mathop {\mathrm {Gms}}\nolimits}

\def\0{\mathbf 0}

\def\ov{\overline}
\def\un{\underline}
\def\wh{\widehat}
\def\wt{\widetilde}
\def\pol{\twoheadrightarrow}

\renewcommand{\rk}{\mathop {\mathrm {rk}}\nolimits}
\renewcommand{\Aut}{\mathop {\mathrm {Aut}}\nolimits}
\renewcommand{\Re}{\mathop {\mathrm {Re}}\nolimits}
\renewcommand{\Im}{\mathop {\mathrm {Im}}\nolimits}
\newcommand{\sgn}{\mathop {\mathrm {sgn}}\nolimits}

\def\bfa{\mathbf a}
\def\bfb{\mathbf b}
\def\bfc{\mathbf c}
\def\bfd{\mathbf d}
\def\bfe{\mathbf e}
\def\bff{\mathbf f}
\def\bfg{\mathbf g}
\def\bfh{\mathbf h}
\def\bfi{\mathbf i}
\def\bfj{\mathbf j}
\def\bfk{\mathbf k}
\def\bfl{\mathbf l}
\def\bfm{\mathbf m}
\def\bfn{\mathbf n}
\def\bfo{\mathbf o}
\def\bfp{\mathbf p}
\def\bfq{\mathbf q}
\def\bfr{\mathbf r}
\def\bfs{\mathbf s}
\def\bft{\mathbf t}
\def\bfu{\mathbf u}
\def\bfv{\mathbf v}
\def\bfw{\mathbf w}
\def\bfx{\mathbf x}
\def\bfy{\mathbf y}
\def\bfz{\mathbf z}

\def\bfA{\mathbf A}
\def\bfB{\mathbf B}
\def\bfC{\mathbf C}
\def\bfD{\mathbf D}
\def\bfE{\mathbf E}
\def\bfF{\mathbf F}
\def\bfG{\mathbf G}
\def\bfH{\mathbf H}
\def\bfI{\mathbf I}
\def\bfJ{\mathbf J}
\def\bfK{\mathbf K}
\def\bfL{\mathbf L}
\def\bfM{\mathbf M}
\def\bfN{\mathbf N}
\def\bfO{\mathbf O}
\def\bfP{\mathbf P}
\def\bfQ{\mathbf Q}
\def\bfR{\mathbf R}
\def\bfS{\mathbf S}
\def\bfT{\mathbf T}
\def\bfU{\mathbf U}
\def\bfV{\mathbf V}
\def\bfW{\mathbf W}
\def\bfX{\mathbf X}
\def\bfY{\mathbf Y}
\def\bfZ{\mathbf Z}

\def\frD{\mathfrak D}
\def\frQ{\mathfrak Q}
\def\frS{\mathfrak S}
\def\frT{\mathfrak T}
\def\frL{\mathfrak L}
\def\frG{\mathfrak G}
\def\frb{\mathfrak b}
\def\frg{\mathfrak g}
\def\frh{\mathfrak h}
\def\frf{\mathfrak f}
\def\frk{\mathfrak k}
\def\frl{\mathfrak l}
\def\frm{\mathfrak m}
\def\frn{\mathfrak n}
\def\fro{\mathfrak o}
\def\frp{\mathfrak p}
\def\frq{\mathfrak q}
\def\frr{\mathfrak r}
\def\frs{\mathfrak s}
\def\frt{\mathfrak t}
\def\fru{\mathfrak u}
\def\frv{\mathfrak v}
\def\frw{\mathfrak w}
\def\frx{\mathfrak x}
\def\fry{\mathfrak y}
\def\frz{\mathfrak z}

\def\bfw{\mathbf w}
%%% END MATHBF
%%%%%%%%%%%%%%%%%%%%%%%%%%%%%%%
%%%%%%%%%%%%%%%%%%%%%%%%%%%%%%%%%
%%% BEGIN MATHBB

\def\R {{\mathbb R }}
 \def\C {{\mathbb C }}
  \def\Z{{\mathbb Z}}
  \def\H{{\mathbb H}}
\def\K{{\mathbb K}}
\def\N{{\mathbb N}}
\def\Q{{\mathbb Q}}
\def\A{{\mathbb A}}
\def\O {{\mathbb O }}

\def\T{\mathbb T}
\def\P{\mathbb P}
\def\SS{\mathbb S}

\def\G{\mathbb G}

\def\cD{\EuScript D}
\def\cL{\EuScript L}
\def\cK{\EuScript K}
\def\cM{\EuScript M}
\def\cN{\EuScript N}
\def\cP{\EuScript P}
\def\cT{\EuScript T}
\def\cQ{\EuScript Q}
\def\cR{\EuScript R}
\def\cW{\EuScript W}
\def\cY{\EuScript Y}
\def\cF{\EuScript F}
\def\cG{\EuScript G}
\def\cZ{\EuScript Z}
\def\cI{\EuScript I}
\def\cJ{\EuScript J}
\def\cB{\EuScript B}
\def\cA{\EuScript A}
\def\cE{\EuScript E}
\def\cC{\EuScript C}
\def\cS{\EuScript S}

\def\bbA{\mathbb A}
\def\bbB{\mathbb B}
\def\bbD{\mathbb D}
\def\bbE{\mathbb E}
\def\bbF{\mathbb F}
\def\bbG{\mathbb G}
\def\bbI{\mathbb I}
\def\bbJ{\mathbb J}
\def\bbL{\mathbb L}
\def\bbM{\mathbb M}
\def\bbN{\mathbb N}
\def\bbO{\mathbb O}
\def\bbP{\mathbb P}
\def\bbQ{\mathbb Q}
\def\bbS{\mathbb S}
\def\bbT{\mathbb T}
\def\bbU{\mathbb U}
\def\bbV{\mathbb V}
\def\bbW{\mathbb W}
\def\bbX{\mathbb X}
\def\bbY{\mathbb Y}

\def\kappa{\varkappa}
\def\epsilon{\varepsilon}
\def\phi{\varphi}
\def\le{\leqslant}
\def\ge{\geqslant}

\def\B{\mathrm B}

\def\la{\langle}
\def\ra{\rangle}
\def\tri{\triangleright}

\def\lambdA{{\boldsymbol{\lambda}}}
\def\alphA{{\boldsymbol{\alpha}}}
\def\betA{{\boldsymbol{\beta}}}
\def\mU{{\boldsymbol{\mu}}}

\def\const{\mathrm{const}}
\def\rem{\mathrm{rem}}
\def\even{\mathrm{even}}
\def\SO{\mathrm{SO}}
\def\OO{\mathrm{O}}
\def\SL{\mathrm{SL}}
\def\PSL{\mathrm{PSL}}
\def\cont{\mathrm{cont}}

\def\U{\operatorname{U}}
\def\GL{\operatorname{GL}}
\def\Mat{\operatorname{Mat}}
\def\End{\operatorname{End}}
\def\Mor{\operatorname{Mor}}
\def\Aut{\operatorname{Aut}}
\def\inv{\operatorname{inv}}
\def\red{\operatorname{red}}
\def\Ind{\operatorname{Ind}}
\def\Fr{\operatorname{Fr}}
\def\dom{\operatorname{dom}}
\def\im{\operatorname{im}}
\def\md{\operatorname{mod\,}}
\def\indef{\operatorname{indef}}
\def\Gr{\operatorname{Gr}}
\def\Pol{\operatorname{Pol}}
\def\Cut{\operatorname{Cut}}
\def\Add{\operatorname{Add}}
\def\ord{\operatorname{ord}}
\def\Replace{\operatorname{Replace}}

\def\arr{\rightrightarrows}
\def\bs{\backslash}

\def\cH{\EuScript{H}}
\def\cO{\EuScript{O}}
\def\cQ{\EuScript{Q}}
\def\cL{\EuScript{L}}
\def\cX{\EuScript{X}}

\def\Di{\Diamond}
\def\di{\diamond}

\def\fin{\mathrm{fin}}
\def\ThetA{\boldsymbol {\Theta}}

\def\0{\boldsymbol{0}}

\def\F{\,{\vphantom{F}}_2F_1}
\def\FF{\,{\vphantom{F}}_3F_2}
\def\H{\,\vphantom{H}^{\phantom{\star}}_2 H_2^\star}
\def\HH{\,\vphantom{H}^{\phantom{\star}}_3 H_3^\star}
\def\Ho{\,\vphantom{H}_2 H_2}

\def\disc{\mathrm{disc}}
\def\cont{\mathrm{cont}}

\def\Kop{\text{\bf\Koppa}}
%\def\Kop{{\bf Q}}

%\Asterisk

%\SixFlowerOpenCenter

%\SixFlowerPetalRemoved

%\def\bigs{\text{\Huge $ *^{\vphantom |}\vphantom{|}$}}

\def\bigs{\AsteriskCenterOpen}
%\def\bigs{$\bigstar$}
%\def\Asterisk{$\bigstar$}

%\def\bigs{\textasteriskcentered}

%\newcommand{\SVER}{\mathop{\Asterisk}\nolimits}

%\Koppa $S$ $\sum$

%$\Kop^3$ $S^3$ $\Kop_3$, $S_3$

%$a$ \XSolid $B$

%$$\epsfbox{chips-add.1}
%$$

%$$\epsfbox{ribbon.1}
%$$

%$$\epsfbox{ribbon.2}
%$$

%$$\epsfbox{ribbon.3}
%$$

%$$\epsfbox{ribbon.16}
%$$

\begin{center}
	\bf \Large
	Virtual permutations and polymorhisms
	
	\bigskip
	
	\sc Yury A.Neretin%
	\footnote{Supported by the grant FWF, Project  P31591.}
	
	\end{center}

%$$\epsfbox{chips-add.5}
%$$

	{\small There is a natural map from a symmetric group $S_n$ to a smaller symmetric group $S_{n-1}$,  we write a decomposition of a permutation into a product of disjoint
	cycles and	 remove  the element $n$ from this expression. For this reason there exists the inverse limit $\mathfrak{S}$ of sets $S_n$. We equip $S_n$ with the uniform distribution (or more generally
		with an Ewens distribution) and get a structure of
		 a measure space on $\mathfrak{S}$ (it is called 'virtual permutations' or 'Chinese restaurant process'), a double $S_\infty\times S_\infty $ of an infinite symmetric group acts  on $\mathfrak{S}$ by left and right 'multiplications'. We discuss the closure of $S_\infty\times S_\infty $ in the semigroup of polymorphisms (spreading maps with spreaded Radon--Nikodym derivatives) of $\mathfrak{S}$. We get formulas for some polymorphisms, in particular for the center of the closure. Expressions  are sums of multiple convolutions of  Dirichlet distributions,  summation sets are certain collections of dessins d'enfant.}
		 
		 	\section{Introduction.  Virtual permutations, \\ chips, and polymorphisms}
		 
		 \COUNTERS

	%$$\epsfbox{chips-add.3}
	%$$
	
	%$$\epsfbox{ribbon.12}
	%$$

	The topic of the paper is formulas for the action of the Olshanski semigroups of chips on the virtual 
	permutations by polymorphims, in this section we explain  these notions  and formulate the problem.
	
	\sm 
	
	%$$\epsfbox{ribbon.15}
	%$$
	
	%$$\epsfbox{chips-add.4} \qquad \epsfbox{ribbon.14}$$
	
	%$${\rm a)}\epsfbox{ribbon.8}\qquad\qquad
%{\rm b)} \epsfbox{ribbon.9}  
%$$

%$${\rm c)}\epsfbox{ribbon.10}\qquad\qquad {\rm d)}  \epsfbox{ribbon.11} 
%$$
%$$
%\epsfbox{ribbon.12} 
%$$

	{\bf\punct Some notation.}
	 Denote by  $\Z_+$ (resp., $\R_+$) the set of nonnegative
	integers (resp. reals), by 
	 $\R^\circ$ (resp. $\R_{>0}$)
	 the  multiplicative group (resp. the set) of positive reals. For a set $A$ we denote by $\#A$ the number of its elements. By $A\sqcup B$, $\sqcup_j A_j$ we denote  disjoint unions of sets.
	 
	 \sm 
	 
	 We denote a list  $(v_1,v_2,\dots)$ by $\{v_i\}$ or by $\{v_i\}_i$
	 (we apply this notation for both ordered or nonordered lists). For finite nonordered
	 lists $\{k_j\}$ of positive integers 
	 we use a  notation 
	 \begin{equation}
	 \iota_m[\{k_j\}]:=\text{number of entries of $m$ to a list $\{k_j\}$.}
	 \label{eq:iota}
	 \end{equation}

	 By $\delta_X[[a]]$ we denote the unit atomic measure (delta-function) at a point $a$ of a space $X$.
	 We denote product-measures by $\mu_1\dot\times \mu_2$ (to not  be confused with the symbol $\times$ of a multiplication).
	 
	 \sm

	Denote by $S_n$ the symmetric group, which is considered as the group of permutations
	of the set 
	$$I_n:=\{1,2,\dots,n\},$$
	 by $S_\infty$ the group of all permutations of $\N$
	with finite supports, 
	$$S_\infty=\lim\limits_{\longrightarrow} S_n=\cup_{n=1}^\infty S_n.$$ 
	This group  is countable and is equipped with the discrete topology.

	By $\ov S_\infty$ we denote the group of all permutations of $\N$.
	This group is continual and is equipped with the unique separable 
	topology compatible with the structure of the group. A sequence $\sigma_j\in \ov S_\infty$
	converges to $\sigma$ if for any $p\in\N$ we have $\sigma_j(p)=\sigma(p)$ for sufficiently
	large $j$. 
	
	\sm 
	
	We can regard  elements of  symmetric groups as  diagrams of the  form
	\begin{equation}
	\epsfbox{chips.1}
	\label{eq:permutation}
	\end{equation}
	We draw a  collection of 'upper circles' enumerated by natural numbers and a collection 
	of 'lower circles' also enumerated by $\N$. If $g$ sends $i$ to $j$, then we connect by an arc
	$i$-th upper circle with $j$-th lower circle. 
	
	A product of permutations corresponds to  a gluing
	of diagrams 
	$$
	\epsfbox{chips.2}
	$$
	
	\vspace{22pt}
	
	{\bf\large A. Virtual permutations and the bisymmetric group}
	
	\vspace{22pt}
	
	{\bf\punct The map $\Upsilon$.%
	\label{ss:virtual}}
Define a canonical map 
$$\Upsilon^n_{n-1}:S_n\to S_{n-1}$$
in the following way. Let $\sigma\in S_n$.

\sm 

1) If $\sigma (n)=n$, we set $\Upsilon^n_{n-1}\sigma(j)=\sigma(j)$ for all $j< n$.

\sm 

2) Let $\sigma(n)=\alpha\ne n$. Then $\sigma^{-1}(n)=\beta$ also $\ne n$.
Then  we assume $\Upsilon^n_{n-1}\sigma(\beta)=\alpha$.
For $j\ne \beta$ we set $\Upsilon^n_{n-1}\sigma(j)=\sigma(j)$.  

\sm

{\sc Example.} For instance,
\begin{align*}
\Upsilon^5_4\begin{pmatrix}
1&2&3&4&\mathbf 5\\3&2&4&1&\mathbf 5
\end{pmatrix}
&=\begin{pmatrix}
1&2&3&4\\3&2&4&1
\end{pmatrix};
\\ \qquad
\Upsilon^5_4\begin{pmatrix}
1&2&\boxed{3}&4&\mathbf 5\\3&2&\mathbf 5&1&\boxed{4}
\end{pmatrix}&=
\begin{pmatrix}
1&2&\boxed{3}&4\\3&2&\boxed{4}&1
\end{pmatrix}.
\end{align*}
On the language of diagrams we have:
$$
\epsfbox{chips.6}
$$
Namely, we add an arc connecting the top '5' and the bottom '5'. We get a compound arc 
$3^{\text{top}}$--$5^{\text{bottom}}$--$5^{\text{top}}$--$4^{\text{bottom}}$ and consider it as an arc 
$3^{\text{top}}$--$4^{\text{bottom}}$ (here and below we consider arcs with fixed ends  up to  isotopies fixing ends).
Similarly,
$$
\epsfbox{chips.9}
$$
We add an arc $5^{\text{ top}}$--$5^{\text{ bottom}}$, get a cycle, and remove it.
		\hfill $\boxtimes$

\sm 

The map $\Upsilon^n_{n-1}$ is $S_{n-1}\times S_{n-1}$-equivariant in the following sense:
\begin{equation}
\Upsilon^n_{n-1} (h_1^{-1}gh_2)=h_1^{-1}\Upsilon^n_{n-1} (g)h_2,\qquad \text{where $g\in S_n$, $h_1$, $h_2\in S_{n-1}$.}
\label{eq:equi}
\end{equation}

	{\sc Another description of the map $\Upsilon^n_{n-1}$.}
	   Decomposing $g\in S_n$
	 as a product of disjoint cycles,
	 \begin{equation}
	 g=(k^1_1\, k^1_2\, \dots)\,(k^2_1\, k^2_2\,\dots)\,\dots
	 \label{eq:decomposition}
	 \end{equation}
	 and removing $n$ from this expression, we get an element 
	 $\Upsilon^n_{n-1}(g)$ defined as product of disjoint cycles. 
	 
	 \sm 
	 
	 More generally, for $n\ge m$ we define a canonical map $\Upsilon^n_m:S_n\to S_m$,
	 $$
	   \Upsilon^n_m(g):=\Upsilon^m_{m-1}\dots  \Upsilon^{n-1}_{n-2}\Upsilon^n_{n-1}(g).
	 $$
	 Equivalently, we remove $m+1$, $m+2$, \dots, $n$ from  expression (\ref{eq:decomposition}) and get an element of $S_m$.
	 
	 \sm

	{\bf \punct The inverse limit of the spaces $S_n$.}
	Thus we have the following chain of surjective maps
	\begin{equation}
	\dots \stackrel{\Upsilon^n_{n-1}}\longleftarrow S_n
	\stackrel{\Upsilon^{n+1}_{n}} \longleftarrow S_{n+1}
	\stackrel{\Upsilon^{n+2}_{n+1}} \longleftarrow S_{n+2}
	\stackrel{\Upsilon^{n+3}_{n+1}} \longleftarrow \dots
	.
	\label{eq:projective}
	\end{equation}
	So we can consider the inverse limit 	of sets $S_n$
	$$
	\frS=\lim\limits_{\longleftarrow} S_n.
	$$
 By the definition, $\frS$ consists of sequences
	$$
	\wt\sigma=(\sigma_1, \sigma_2, \sigma_3,\dots), 
	\qquad
	\text{where $\sigma_k\in S_k$ and $\Upsilon^k_{k-1} \sigma_k=\sigma_{k-1}$.}
	$$
	In particular, we have a canonical map
	$\Upsilon^\infty_n:\frS\to S_n$ defined by $\Upsilon^\infty_n \wt\sigma=\sigma_n$.

	\sm 
	
	A point of the space $\frS$ can be described in the following way. We fix a  partition
	of $\N$ into a disjoint union of subsets ({\it pre-tables})
	\begin{equation}
	\N=M_1\sqcup M_2\sqcup\dots
	.
	\label{eq:N-cycles}
	\end{equation}
	and fix a cyclic ordering in each subset.
	A map $\Upsilon^\infty_n$ is the removing all elements $>n$. Then we
	get the set $I_n=\{1,2,\dots,n\}$ splitted into a disjoint union of cyclically ordered subsets, i.e., an element of $S_n$.

	\sm
	
	The space $\frS$ is not an inverse limit in the category of groups, since generally 
	$$
	\Upsilon^n_{n-1}(g_1 g_2)\ne \Upsilon^n_{n-1}(g_1)\Upsilon^n_{n-1}(g_2).
	$$
	In particular, $\frS$ does not have a group structure.
	
	 However, for $n>m$
	by \eqref{eq:equi},
		$$
	\Upsilon^n_{m}(h_1^{-1}g h_2)=h_1^{-1} \Upsilon^n_{m}(g) h_2,
	\qquad 
	\text{where $g\in S_n$, $h_1$, $h_2\in S_{m}$.}
	$$

	Therefore $S_m\times S_m$ acts on $\frS$.
	 Since $m$ is arbitrary, we get an action of the double $S_\infty\times S_\infty$
	on $\frS$.
	% Namely,
	%\begin{multline*}
	%h_1\wt \sigma h_2^{-1}:=
	%\\
	%:=\bigl(\Upsilon^n_1(h_1^{-1}\sigma_n h_2), \dots, %\Upsilon^n_{n-1}(h_1^{-1}\sigma_n h_2),
	%h_1^{-1}\sigma_n h_2, h_1^{-1}\sigma_{n+1} h_2, %h_1^{-1}\sigma_{n+2} h_2,\dots\bigr)
	%\end{multline*}
	
	 It is more natural to consider 
	such inverse limit  as a measure space.
	
	\sm 
	 
	 {\bf\punct Virtual permutations.%
	 \label{ss:virt}}
	 	For $g\in S_n$ denote by $[g]$ the number of its disjoint cycles. 
	 	Fix $z>0$.
	 {\it The Ewens distribution} is a probabilistic measure $\mu^z_n$  on $S_n$ defined by 
	 $$
	 \mu^z_n(g)=\frac{z^{[g]}}{z(z+1)\dots(z+n-1)},
	 $$
	 for $z=1$ we get a uniform distribution on $S_n$.
	 
	 It easy to show that the pushforward of the Ewens distribution $\mu^z_n$
	 is the Ewens distribution $\mu^z_{n-1}$ (i.~e., for any subset $K\subset S_{n-1}$
	 the measure of its preimage $(\Upsilon^n_{n-1})^{-1}(K)$
	 coincides with the measure of $K$). We  define the space $\frS^z$
	 of {\it virtual permutations} as the  inverse limit
	 \begin{equation}
	 (\frS^z,\mu^z):=\lim_{\longleftarrow} (S_n,\mu^z_n)
	 \label{eq:inverse}
	 \end{equation}
	in the category of measure spaces. Namely, we equip $\frS$ with
	 the measure $\mu^z$   defined from the condition:
	for any $n$ and $g\in S_n$ we set $\mu^z((\Upsilon^\infty_n)^{-1}g)=\mu_n^z(g)$. By the Carath\'eodory extension theorem we get a well-defined measure on $\frS$.
%	It is more or less clear that this space as a measure space is equivalent to a segment $[0,1]$ equipped with the Lebesgue measure.

\sm 

	{\sc Remark.} Apparently, a construction of $\frS^z$ was firstly appeared in Aldous \cite{Ald} (1985), 11.19,
		 with a reference to J.~Pitman (see, also, Pitman \cite{Pit}). The action of the group $S_\infty\times S_\infty$
	on $\frS$ appeared in Kerov, Olshanski, Vershik \cite{KOV0}.
	Notice that for a fixed $z$ for each $n$ we have a distribution of lengths of cycles on $S_n$, i.e., a distribution
	on the set of partitions of $n$. So we have a Markov process of growth
	of such partitions, this process arises to the work of Ewens \cite{Ewe}, 1972, on population genetics.
	\hfill  $\boxtimes$

	\sm

	{\bf \punct A more convenient description of $\frS^z$.}
	It is done in two steps.
	
	\sm 
	
	{\sc Step 1.} We present the limit distribution of lengths
	of cycles, see Kingman \cite{Kingman}.
	Consider the measure $\nu$
	on the half-line $x>0$ defined by
	$$
	d\nu^z(x)=z\, x^{-1}e^{-x}\, dx.
	$$
	Consider the Poisson measure $\pi^z$  on the set of (non-ordered) countable subsets (see, e.g., \cite{Kingman},
	\cite{Ner-book}, Sect. X.4)
	$\{s_1,s_2,\dots \}$ in $\R$ defined by the measure%
	\footnote{For a measurable subset $A\subset \R_{>0}$ denote by $\Omega (A,k)$ the set ('event') of all configurations whose
	intersections with $A$ have precisely $k$ point. We assume that a propbabibility of $\Omega(A,k)$ is $\exp\{-\nu^z(A)\}\, \nu^z(A)^k/k!$
	and for pairwise disjoint sets $A_j$ and arbitrary $k_j$ the events $\Omega(A_j,k_j)$ are independent.} 
	$\nu^z$. It can be easily shown
	 that the series $\sum s_j$ converges a.s.  
	We define a collection $\{\ell_\omega\}$ as $\ell_\omega=s_\omega/\sum s_j$ and get a measure ({\it the Poisson--Dirichlet distribution}) on the set of unordered collections $\{\ell_\omega\}$ such that $\sum \ell_\omega=1$.

	For a given $\{\ell_\omega\}$ consider a collection $U_1$, $U_2$, $\dots$ of oriented circles of lengths $\ell_1$, $\ell_2$, \dots . We  imagine them as circles on plane oriented clockwise and call them by {\it tables}, sets $\sqcup U_\omega$ we call by  {\it restaurants} ('Chinese restaurants').
	Denote the space of all restaurants equipped with the Poisson--Dirichlet distribution by $\frT^z$.
	
	We can regard points of $\frT^z$ as non-ordered collections of tables $\{U_\omega\}$ or as non-ordered collections of lengths $\{\ell_\omega \}$.
	
	\sm 
	
	{\sc Remark.} Clearly, the numbers $\ell_\omega$ are pairwise distinct a.s.
	\hfill $\boxtimes$
	
	\sm
	
	{\sc Step 2.} We chose a random sequence
	in $\sqcup_\omega U_\omega$ such that each element is distributed uniformly with respect
	to the Lebesgue measure  on $\sqcup_\omega U_\omega$.
	 We denote points of a sequence (we call them  {\it guests})
	by $\ov 1$, $\ov 2$, $\ov 3$, \dots. Configurations of guests on  tables are defined upto  rotations of  tables.
	A {\it virtual permutation} or an {\it occupied restaurant} is a collection $\{\ell_\omega\}$
	distributed according the Poisson--Dirichlet law, the corresponding restaurant $\sqcup U_\omega$, and
	a random sequence $\ov 1, \ov 2, \dots\in \sqcup_\omega U_\omega$.
	We denote the measure space of all virtual permutations
	by $\frS^z$, denote by $\mu^z$ the measure on $\frS^z$.
	
	We denote {\it occupied tables} by $\ov U_\omega$.

	\sm 
	
	{\sc Remark.} Clearly, for almost all points of $\frS^z$  a sequence of guests $\ov 1$, $\ov2$, \dots is dense in $\sqcup_\omega U_\omega$
	and its elements are pairwise distinct. 
	\hfill $\boxtimes$
	
	\sm 
	
{\sc Remark.} The set  of all  tables  of a given restaurant can be identified with  $\N$, 
	we can order sets of tables according a length $\ell_\omega$ or according  a minimal number of a guest sitting at the table.
	We prefer
	to think that a set of tables of a given restaurant is an    abstract countable set and do not identify such sets for different restaurants. 
	\hfill $\boxtimes$

	\sm 
	
	There is a canonical map $\Upsilon^\infty_n:\frS^z\to S_n$. Namely, let us 'forget' all guest $\ov{n+1}$, $\ov{n+2}$, \dots $\in\sqcup  U_\omega$. After this, each table $U_\omega$ contains a finite collection of guests 
	$\ov{i_\omega^1}$, $\ov{i_\omega^2}$, \dots, $\ov{i_\omega^{k_\omega}}$
	seating at the table. So we get a partition of $\{1,\dots,n\}$ into a disjoint union of cyclically ordered
	sets.  Thus we come to a permutation decomposed into a product of disjoint  cycles.
	
	The pushforward  of the measure $\mu^z$ to $S_n$ is the Ewens measure
	$\mu^z_n$. We also have $\Upsilon^n_{n-1}\circ\Upsilon^\infty_n=\Upsilon^\infty_{n-1}$.
	
	\sm
	
	It  remains to explain how pre-tables defined in Subsect. \ref{ss:virt} generate tables. Let $\wt\sigma$ be a point of the 
	inverse limit (\ref{eq:inverse}).  Consider two points $i$, $j$
	on one pre-table.
	For each $N\ge i,j$ consider the cycle of 
	$\Upsilon^\infty_N \wt \sigma$ containing $i$, $j$, say $(\dots\,i\, m_1\, m_2\,\dots\, m_p\, j\, \dots)$. Denote by $p_{ij}(N)$ the number of elements between
	$i$ and $j$.
		Then (see \cite{Tsi}, \cite{KT}) for almost all points of $\frS^z$ for all 
		pairs $\ov i$, $\ov j$ lying on one pre-table we define
	\begin{equation}
	\text{length of arc $[\ov i, \ov j]$}:=
	\lim_{N\to\infty} \frac{p_{ij}(N)}{N}.
	\label{eq:limit}
	\end{equation}
	 The limit exists a.s. for all $i$, $j$. Then we get an occupied restasurant.

	\sm 

{\bf \punct The action of $S_\infty\times S_\infty$ on $\frS$.%
\label{ss:action-group}}	
	The space $\frS^z$ is not a group, but the group $S_\infty\times S_\infty$
	acts on $\frS$ by 'left and right multiplications'.  Let $\fru\in \frS^z$, 
	let $\{\ov U_\alpha\}$ be its occupied  tables. Set $u=\Upsilon^\infty_n \fru$.
	For each $i\in I_n$ we denote by $\ov U[i]$ the table containing the guest $\ov i$,
	by $\ov U_+[i]$ the arc $[\ov i, \ov{u(i)}]$ of $\ov U[i]$, by  $\ov U_-[i]$ the arc $[\ov{u^{-1}(i)}, \ov i]$.

	Fix $g\in S_n$.  Cutting tables at points (guests) $\ov 1$, \dots, $\ov n$, we get a finite collection of segments $\{\ov U_+[i]\}_{i\in I_n}$  and a countable number of non-cutted tables. The collection $\{\ov U_-[j]\}_{j\in I(n)}$ coincides with $\{\ov U_+[i]\}_{i\in I_n}$ up to a reordering.
	
	To obtain $\fru g$, for each $i\in I_n$ we glue the segments $\ov U_-[j]$ and $\ov U_+[g(j)]$
	identifying points $\ov j\in \ov U_-[j]$ and $\ov{g(j)}\in\ov U_+[g(j)]$ and putting the guest $\ov j$
	to the point of gluing. In this way, we get a family of new tables enumerated by disjoint cycles of $ug\in S_n$, and add the collection of non-cutted tables.
	
    To obtain $g\fru$ we repeat the same steps, but put the guest $\ov{g(j)}$ to the point of gluing. Notice that $g\fru =g(\fru g) g^{-1}$, so the restorants
    of $\fru g$ and $g\fru$ have the same  tables and the same positions of all  guests with numbers $>n$, the guests $\ov 1$, \dots, $\ov n$ differ by a permutation $g$.

 %Let us cut the cycles $L_\alpha$ of $u$ at points $\ov1$, \dots, $\ov n$
%	and get a finite collection of directed segments $[{i_\alpha^m},i_\alpha^{m+1}]$.
%	The set of all origins is precisely $\{\ov1,\dots,\ov n \}$, the set of all ends is 
%	also $\{\ov1,\dots,\ov n \}$. Also we have a countable collection of survived cycles.

%---	To obtain $hu$, for each $p\le n$ we found a unique $p$ at the ends of segments, 
%	say $[\ov{i_\beta^1},\dots,\ov{i_\beta^{k_\beta-1}},\ov p]$ and found a unique $h(p)$
%	at the origins of segments. Next, we identify points $p$ at the ends with corresponding
%	points $h(p)$ on the origins, and put the labels $\ov p$ on such points. After this we get
%	a countable collection of cycles with labels $\ov 1$, $\ov 2$,\dots

%--- To obtain $uh$, for each $q\le n$	we found $h (q)$ among origins of segments and $q$ on the origins. We identify each $h(q)$ among the beginnings with each $q$ on origins. Next, we identify
%each $h(q)$ at a beginning with the corresponding $q$ at the end and put the label $\ov q$ on each point of gluing.

\sm 

{\sc Remark.} These rules are simply a rephrasing of a description of  multiplication of a permutation defined as map and a permutation
defined as a product of cycles.
\hfill $\boxtimes$

\sm

{\bf \punct  Radon--Nykodim derivatives.}
Thus, we get an action of the group $S_\infty\times S_\infty$ on $\frS^z$ given by
$\fru\mapsto h_1^{-1}\fru h_2$.

Clearly, on the  level of $S_n$ 
we have
$$
\mu^z_n(h_1^{-1}u h_2)= z^{[h^{-1}_1gh_2]- [g]} \mu_n^z(u).
$$
For $\fru:=\frS^z$ the operation of cutting and gluing described above
involves only a finite number of tables, so the number
\begin{multline*}
\gamma(h_1,h_2;\fru):=\bigl\{\text{number of tables of $h_1^{-1}\fru h_2$}\bigr\} - \bigl\{\text{number of tables of $\fru$}\bigr\}
:=\\
:=[h^{-1}\Upsilon^\infty_N(g) h_2]- [\Upsilon^\infty_N(g)]\qquad \text{for large $N$}
%:=\bigl\{\text{number of cycles of $h_1^{-1}u h_2$}\bigr\} - \bigl\{\text{number of cycles of $u$}\bigr\}
\end{multline*}
is well-defined. The Radon-Nikodym derivative%
\footnote{See a formal definition below in Subsect. \ref{ss:poli}, see also, e.g., \cite{Bog}, 9.12.}
of the transformation
	$\fru\mapsto h_1^{-1}\fru h_2$
is $z^{\gamma(h_1,h_2;\fru)}$.

\sm

{\bf \punct The inversion of virtual permutations.}
Clearly,
$$
\Upsilon^n_{n-1}(g^{-1})=g^{-1}.
$$
So the inversion $g\mapsto g^{-1}$ defines a measure preserving map $\frS^z\to\frS^z$, we denote it by $\fru\mapsto \fru^{-1}$. This is simply
a changing of orientations of all tables from clockwise to counterclockwise. Clearly
$$
(h\fru g)^{-1}=g^{-1}\fru^{-1}h^{-1}
.$$  

{\bf \punct 	The bisymmetric group.} The {\it bisymmetric group} 
$\SS$ is the subgroup in $\ov S_\infty \times \ov S_\infty$ consisting of pairs $(g_1,g_2)$ such that
$g_1g_2^{-1}\in S_\infty$. Denote by $\K\subset \SS$ the diagonal subgroup, i.e., the subgroup
consisting of pairs $(g,g)$, where $g$ ranges in $\ov S_\infty$. So, $\K\simeq \ov S_\infty$.
We define topology on $\SS$ assuming that $\K$ is an open subgroup equipped with the natural topology.
The homogeneous space $\SS/\K$ is a countable space with the discrete topology.

\sm 

{\it The bisymmetric group $\SS$ acts on $\frS^z$}. Namely, the subgroup $\K$ acts by permutations of guests, such permutations are
measure-preserving maps. The subgroup $S_\infty\times S_\infty$
acts as above.

\sm 

{\sc Remark.}
1) The group $\SS$ was introduced by Olshanski  \cite{Olsh-symm}. It is a type $I$ group, a classification
of its irreducible unitary representations is known, see Olshanski \cite{Olsh-symm},  Okounkov \cite{Oku}, see also \cite{Mel}.
The most important family of irreducible  representations are $\K$-spherical representations, i.e., representations having a fixed vector with respect
to $\K$ (it is automaticaly unique). In fact such representations were classified by 
Thoma \cite{Tho}, 1964, see also \cite{VK} and \cite{Olsh-symm}.

\sm 

2) Decompositions of quasiregular representations of $\SS$ in $L^2(\frS^z)$ are known,
see Kerov, Olshanski, Vershik \cite{KOV}, Borodin, Olshanski
 \cite{BO},  \cite{BO-book}. One of informal aims 
of the present work is a search of additional possibilities for the  analysis on the space     
 $\frS^z$.

	\vspace{22pt}

{\bf\large B. The bisymmetric group and its train.}

\vspace{22pt}

{\bf \punct Multiplication of double cosets.} See Olshanski \cite{Olsh-symm}, see also \cite{Ner-umn}.
	For $\alpha\in\Z_+$ we denote by $\K_\alpha$ the subgroup  in $\K$ consisting of $(g,g)$ such that
$g\in \ov S_\infty$ fixes $1$, \dots, $\alpha\in \N$. We set $\K_0:=\K$. Denote by
$$K_\alpha:=\K_\alpha\cap \bigl(S_\infty\times S_\infty \bigr)$$
the corresponding subgroup of finitely supported permutations.

 For any $\alpha$, $\beta$ we consider the double coset space 
$$
\cS(\alpha,\beta)=\K_\alpha\backslash \SS/\K_\beta.
$$
It is easy to see that each double coset $\K_\alpha \cdot g\cdot \K_\beta$
has a finitary representative $g'\in S_\infty\times S_\infty$.
Moreover,
\begin{equation}
\K_\alpha\backslash \SS/\K_\beta\simeq K_\alpha\backslash  \bigl(S_\infty\times S_\infty \bigr)/K_\beta,
\label{eq:two-cosets}
\end{equation}
 all such quotient spaces are countable and the quotient topologies are  discrete.

\begin{figure}
	$$
	\epsfbox{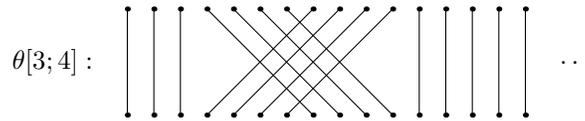}
	$$
	\caption{To the definition of a product of double cosets.
%We draw an element of $S_\infty$ as two rows of black circles and a collection
%of arcs connected circles of the top row and circles in the bottom row.	
}
	\label{fig:razdviganie}
\end{figure}

For any $\alpha$, $\beta$, $\gamma\in \Z_+$,
there exists a natural multiplication %$\frg\circ\frh$
$$
\cS(\alpha,\beta)\times \cS(\beta,\gamma)\to \cS(\alpha,\gamma)
$$	
defined in the following way. For each $\beta\in \Z_+$
we consider the following sequence $\theta^\beta[j]\in \K_\beta$ 
defined by
$$
\theta^\beta[j](k)=
\begin{cases} 
\alpha &\qquad\text{if $k\le \beta$;}
\\
\alpha+j&\qquad\text{if $\beta <k\le \beta +j$;}
\\
\alpha-j&\qquad\text{if $\beta+j<k\le \beta+2j$;}
\\
\alpha&\qquad \text{if $k>
	\beta+2j$.}
\end{cases},
$$
see Fig.~\ref{fig:razdviganie}.

Now we take double cosets $\frg=\K_\alpha g \K_\beta$, $\frh=\K_\beta h \K_\gamma$,
without loss of generality we can think that $g$, $h\in S_\infty\times S_\infty$.
Consider the sequence 
$$
\K_\alpha \cdot g\, \theta^\beta[j]\, h \cdot \K_\gamma\,\,\in \K_\alpha \backslash \SS/\K_\beta.
$$
It is easy to show that this sequence is eventually constant. We define a product
$\frg \circ \frh\in \K_\alpha\backslash \SS/\K_\gamma$ as the limit (i.e., the stable value) of this sequence. It can be shown that the result
does not depend on a choice of representatives $g$, $h\in S_\infty\times S_\infty$
and this operation is associative. So we get a category $\cS$
({\it the train of $\SS$}) whose objects
are nonnegative integers and morphisms $\beta\to \alpha$ are double cosets,
$$
\Ob(\cS):=\Z_+,\qquad \Mor_\cS(\beta,\alpha):=\cS(\alpha,\beta)=
\K_\alpha\backslash \SS/\K_\beta.
$$

We also define an involution $\frg\mapsto \frg^*$ on the category $\cS$. Namely, a map $g\mapsto g^{-1}$
determines bijections $\K_\alpha\backslash \SS/\K_\beta\to \K_\beta\backslash \SS/\K_\alpha$. Clearly,
$$
(\frg\circ\frh)^*=\frh^*\circ \frg^*.
$$

\sm 

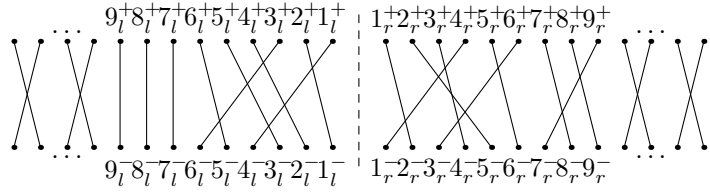
\begin{figure}
	$$
	\epsfbox{chips.3}
	$$
	a) An element of $\SS$. The symbols $l$ and $r$ in subscripts are abbreviations 
of 'left' and 'right'.  Symbols $+$ and $-$ in superscripts correspond to the top and the bottom of the diagram.
 The left and the right parts 
	of the diagram are symmetric one to another except a finite number of arcs.
	
	$$
	
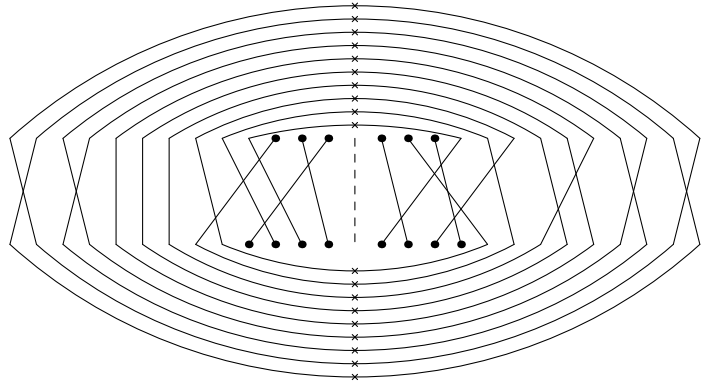
\epsfbox{chips.4}
	$$
	b) On Figure $\beta=3$, $\alpha=4$. We add 'horisontal' arcs with crosses.
	
	$$
	
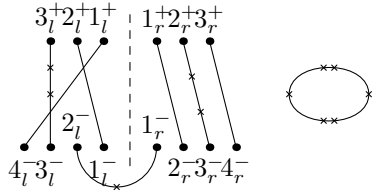
\epsfbox{chips.5}
	$$
	c)    We consider compound arcs on the previous diagram
	up to isotopies with fixed ends. Also, we get a countable family of circles with two crosses and forget  them. 
	
	\caption{The construction of a chip from an element of $\SS$.}
	\label{fig:1}
\end{figure}

{\bf \punct Chips.}
Recall Olshanski's description \cite{Olsh-symm}
 of the product of double cosets.
We represent an element of the group $S_\infty\times S_\infty$ as a diagram 
of the type drawn
on Fig.~\ref{fig:1}.a.
On Fig.~\ref{fig:1}.b, \ref{fig:1}.c we define combinatorial data corresponding to a double coset.
We get a diagram (a {\it chip%
\footnote{We use the term proposed by Kerov \cite{Ker} for elements of the Brauer semigroups.
	% Recall that figures similar to Fig.~\ref{fi:1} appeared in the work \cite{Br}.
	 Recall that the Schur--Weyl duality between
$\GL(N)$ and $S_n$ has a counterpart \cite{Br}, 1936, for orthogonal groups $\OO(N)$ and symplectic groups $\Sp(2N)$, dual objects are certain semigroups
of chips (with horisontal arcs but  without crosses). Brauer's  approach also gives a similar
statement  for tensors over $\GL(N)$ with covariant and contravariant components, in this case diagrams have a separating dashed line (it separates 'co-' and 'contra-') and horizontal arcs corresponding to convolutions. On the other hand, according Wasserman
\cite{Was}
and Olshanski \cite{Olsh-symm},
representations of $\SS$ can be realized in certain tensors of infinite order,  their decomposition is controlled by certain 'dual' compact groups.}})
	 of the form
shown on the Fig.~\ref{fig:1}.c.

It is convenient to think that {\it each arc has
a 'length', which is defined as the half of the number of crosses.}

{\it For figures in our printed text  it is more convenient to draw crosses.}

\sm 

We get diagrams of the following type. On the top we have 
 black circles labeled by  $\beta_l^+$, \dots, $1_l^+$ to the left of the dashed line, and 
   $1_r^+$, \dots, $\beta_r^+$ to the right of the dashed line. On the bottom we
   have black circles labeled by $\alpha_l^-$, \dots, $1_l^-$ on the left
and   $1_r^-$, \dots, $\alpha_r^-$ on the right. Each circle is an end of an arc.
There are 3 following  types of arcs:

\sm

1) arcs $i_r^+\downarrow i_r^-[\phi]$ (resp., $i_l^+\downarrow i_l^-[\phi]$) 
 of integer length $\phi\ge 	0$  from the top to the bottom located 
in the right (resp., left) hand side of the diagram; 
\sm 

2) arcs $i_l^+\smallfrown j_r^+[\psi]$ (resp., $i_l^-\smallsmile j_r^-[\psi]$)
of length $\psi\in (1/2+\Z_+)$
from left to right on the top (resp., bottom) of the diagram;

\sm 

3) circles $\bigcirc [k]$ of integer length $k>0$.

\sm

{\it We remove  all cycles of length} 1
(any diagram obtained by our procedure  contains an infinite number of cycles $\bigcirc [1]$, this collection contains no 
information). So in all cases we get finite objects.

\sm

\begin{figure}
	$$\epsfbox{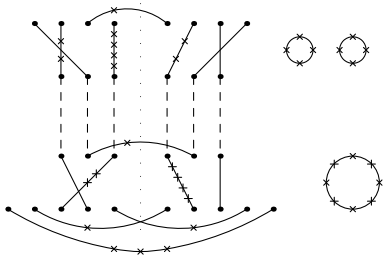}
	$$
	
	\caption{The multiplication of chips. We identify corresponding circles of the bottom of the first diagram and on the top of the second diagram.}
	\label{fig:multiplication}
\end{figure}

It is easy to show that such diagrams are in one-to-one correspondence with elements of $\K_\alpha\backslash \SS/\K_\beta$.
Multiplication of diagrams is a gluing, see 
Fig.~\ref{fig:multiplication}. The involution corresponds to the replacement of the top and of the bottom.

\sm 

{\sc Remark.} An element of the semigroup $\cS(0,0)$ is a collection of cycles $\bigcirc[k_j]$.
The product is a union of such collections. In particular, this semigroup is Abelian.
\hfill $\boxtimes$

{\bf \punct The multiplicativity theorem.}
Let $\rho$ be a unitary representation of the group
$\SS$ in a Hilbert space%
\footnote{We assume that a Hilbert space is separable by definition.} $H$. Denote by $H_\alpha$
the subspace of all $\K_\alpha$-fixed vectors, by $P_\alpha$ the operator of orthogonal projection to $H_\alpha$.
For $g\in \SS$ we define the operator
$$
\wt \rho_{\alpha,\beta}(g): H_\beta\to H_\alpha
$$
 by
$$
\wt \rho_{\alpha,\beta}(g):= P_\alpha\, \rho(g)\Bigr|_{H_\beta}.
$$
It is easy to see that $\wt \rho_{\alpha,\beta}(g)$
depends only on the double coset
$\frg=\K_\alpha g \K_\beta$.
The following {\it multiplicativity theorem} holds (see \cite{Olsh-symm}, \cite{Ner-umn}).

\begin{theorem}
	\label{th:0}
	For any $\alpha$, $\beta$, $\gamma\in \Z_+$
	and $\frg\in \K_\alpha\backslash \SS/\K_\beta$,
	$\frh\in \K_\beta\backslash \SS/\K_\gamma$
	we have
	\begin{equation}
	\wt \rho_{\alpha,\beta}(\frg)
	\,
	\wt\rho_{\beta,\gamma}(\frh)=
	\wt\rho_{\alpha,\gamma}(\frg\circ \frh).
	\label{eq:multiplicativity}
	\end{equation}
\end{theorem} 

\begin{equation}
\wt\rho(\frg^*)=\wt\rho(\frg)^*
\label{eq:invinv}
\end{equation}

So, for any unitary representation of the group 
$\SS$ we get a functor from the category $\cS$ to the category of Hilbert spaces and bounded operators. 
This  reduces investigation of representations of the group $\SS$ to representations of the category $\cS$.

\sm 

In particular, this theorem can be applied
to  quasiregular representations of the bisymmetric group $\SS$ in
$L^2(\frS^z)$. There arises a question: {\it is it possible to 
describe operators corresponding to chips explicitly?}
Below we reformulate this question in terms of polymorphisms.

	\sm
	
{\bf \punct Remark. Infinite chips and weak closures of the group $\SS$ in unitary representations.%
\label{ss:infty-infty}}
Denote by $\cC(H)$ the semigroup of all operators in a Hilbert space $H$ with norm $\le 1$ (we call such operators {\it contractive}). It easy to see  that
$\cC(H)$ is compact and metrizable with respect to the weak operator 
topology, the multiplication in $\cC(H)$
is separately continuous (see, e.g., \cite{Ner-book}, Sect.~I.1).

Let $G$ be a group, $\rho$ be its unitary representation
in a Hilbert space $H$. Consider the subset $\rho(G)$
in the space of all operators in $H$ and consider its closure $\ov{\rho(G)}$ in the semigroup of all contractive operators
in the weak operator topology. Clearly, we  get a compact separately continuous semigroup, 
see \cite{Olsh-semi}, \cite{Ner-book}, Sect.~I.1.  

A description of such semigroup can be a nontrivial problem even for $\Z$, i.e.,
for closure of the set of powers of a given unitary operator $U$.
For a comeagre set of unitary operators this semigroup
is the unit ball in $L^\infty(S^1,\mu)$, where $\mu$ is the spectral measure
of $U$, see \cite{MTs}. 

For reasonable representations of Lie groups (over real and over local fields) the question usually leads to one-point compactifications, 
see \cite{HM}. Other types of locally compact groups and  discrete non-Abelian groups are not well-understood. For infinite-dimensional groups this question leads to handable and unexpected
algebraic structures as it was observed by Olshanski \cite{Olsh-semi}. Let us describe such
compactification for the group $\SS$.

 The sequence $\theta_\beta[j]$ defined above
converges to the projector $P_\beta$ in the weak operator topology, see, e.g., \cite{Ner-book}, Theorem VIII.1.4. 
Therefore the following operators
$$
\wh \rho_{\alpha,\beta}(g)=P_\alpha\,\rho(g)\,P_\beta: \,H\to H
$$
are contained in $\ov{\rho(\SS)}$.
It is easy to see that these operators depend only on double coset $\frg$ containing $g$ and have the following block structure
$$
\wh\rho (\frg)=
\begin{pmatrix}
\wt \rho(\frg)&0\\0&0
\end{pmatrix}:\quad 
H_\beta\oplus H_\beta^\bot\to H_\alpha\oplus H_\alpha^\bot.
$$
Clearly,   we have
$$
	\wh \rho_{\alpha,\beta}(\frg)\,
\wh\rho_{\beta,\gamma}(\frh)=
\wh\rho_{\alpha,\gamma}(\frg\circ \frh),
$$
this identity is a rephrasing of
(\ref{eq:multiplicativity}).
So chips act in the space $H$ itself.

	\sm 
	
More generally, consider diagrams as on Fig.~\ref{fig:1}.c, where circles on the top and the bottom (from the left and from the right) are enumerated by $\N$. We require the following conditions of finiteness of such diagrams:

\sm 

a) the number of cycles $\bigcirc[k]$ is finite
(recall that $k\ge 2$);

\sm 

b) the number of 'vertical' arcs $i_l^+\downarrow j_l^-[m]$,
$p_r^+\downarrow q_r^-[m]$ of lengths $m>0$ is finite;

\sm 

c) the number of horizontal arcs $i_l^+\smallfrown j_r^+[1/2+k]$,
$i_l^-\smallsmile j_r^-[1/2+k]$
of lengths $1/2+k>1/2$ is finite.

\sm 

d)  
a diaram is symmetric with respect to the dashed line upto a finite family of arcs.

\sm 

Denote the semigroup of such chips  by $\cS(\infty,\infty)$. 
The group of its invertible elements is $\SS$. 

\sm

{\sc Remark.} The embedding $\cS(n,n)\to\cS(\infty,\infty)$
is given by adding  arcs $i^+_l\smallfrown i^+_r[1/2]$, $i^-_l\smallsmile i^-_r[1/2]$ for all $i>n$.
\hfill $\boxtimes$

%There is a natural topology on $\Delta(\infty,\infty)$,
%see below \ref{}. The group $\SS$ is dense in $\Delta(\infty,\infty)$.

\sm 

{\sc Remark.} The semigroup $\cS(\infty,\infty)$ has a center. It consists of chips whose arcs have the form
$j_l^+\downarrow j_l^-$,  $j_r^+\downarrow j_r^-$ for all $j\in \N$  and 
$\bigcirc[k_j]$, where $k_j\ge 2$.  The center is isomorphic to $\cS(0,0)$.
\hfill $\boxtimes$

\sm

For a unitary representation $\rho$ of $\SS$ denote by $\Pi_+$ 
the operator of orthogonal projection to the subspace of vectors fixed by the whole group $\SS$, by $\Pi_-$ the operator of projection to
the subspace of vectors $v$ satisfying
$$
\rho(h_1, h_2) v= \sgn(h_1h_2^{-1})v,\qquad\text{where $(h_1, h_2)\in \SS$}.
$$

\begin{proposition}
	\label{pr:closure}
	{\bf a)} Any unitary representation $\rho$ of $\SS$ extends to
	a  representation of the semigroup $\cS(\infty,\infty)$.
	
	\sm 
	
	{\rm b)} The closure $\ov{\rho(\SS)}$ consists of
	the image of $\cS(\infty,\infty)$ and operators
	$\Pi_++\Pi_-$ and $\Pi_+-\Pi_-$. 
\end{proposition}

{\sc Remarks.} a) Representations of the semigroup $\cS(\infty,\infty)$ can be constructed in the following way. We split
$\N$ into a disjoint union of two countable sets $\N=L\sqcup M$. Denote by $\K(\infty)$ the subgroup
of $\K$ consisting of elements fixing all elements of $L$. Then $\cS(\infty,\infty)\simeq \K(\infty)\backslash \SS/\K(\infty)$
and the multiplicativity theorem remains valid with the same proof. The group of invertible elements of  $\cS(\infty,\infty)$ is isomorphic to
 $S(\infty,\infty)$. It can be shown that the  representation of $\SS$ in the subspace $\K(\infty)$-fixed vectors is
equivalent to the initial representation of $\SS$.

\sm

b) The second statement of the theorem is a relatively easy corollary of results of Olshanski \cite{Olsh-symm} and Okounkov \cite{Oku},
we omit its proof.

\sm

c) We can define sets $\cS(\infty,m)$ and $\cS(n,\infty)$ in the same way, and so get a completed category $\ov\cS$,
$$
\Ob(\cS):=\Z_+\sqcup \infty.
$$
 Any unitary representation of $\SS$ admits a unique  extension to a representation of 
this category. \hfill $\boxtimes$

%The statement is a simple corollary of works Olshanski \cite{Olsh-symm} and Okunkov \cite{Oku}, we omit a proof.

\sm

\vspace{22pt}

{\bf C. Polymorphisms}

\vspace{22pt}

It is natural to reformulate the question about operators in $\frS^z$
corresponding to chips in terms of polymorphisms of the space
$\frS^z$.

\sm 

{\bf \punct Bistochastic kernels or measure preserving polymorphisms.%
\label{ss:bist}} Recall that a {\it Lebesgue  space}
is a measure space, which is equivalent to a union of an interval of a line and a finite, countable or empty set of (atomic) points with nonzero measures (see, \cite{Roh}, \cite{Bog}, \S9.4).
A measure $\mu$  is {\it probabilistic} if the measure of the whole space is 1.  A measure is {\it continuous} if it has no atomic points.
Clearly, the spaces $\frS^z$ of virtual permutations are Lebesgue spaces with continuous measures.

Let $(M,\mu)$ be a Lebesgue space with a continuous measure $\mu$. Denote by $\Ams(M)$ the group of bijective a.s. measure preserving transformations of $M$.
% by $\Gms(M)$ 
%the group of transformations $g$ leaving the measure quasiinvariant, 
%i.e., transformations sending $\mu$ to equivalent measures. 
 The group $\Ams(M)$ acts in $L^2(M,\mu)$ by unitary operators
\begin{equation}
T(g)f(m)=f(mg).
\label{eq:T-pres}
\end{equation} 

Let $(M,\mu)$, $(N,\nu)$ be probabilistic Lebesgue measure spaces.
A {\it bistochastic kernel} or a {\it measure preseving polymorphism}%
\footnote{The notion rises to E.~Hopf \cite{Hopf}, see also
\cite{Nev}, \cite{Kre}, \cite{Ver}, such objects are widely 
used in  ergodic theory, see e.g., \cite{King}, \cite{JPR}.
 We use the term 'polymorphism' proposed in \cite{Ver}.}
 $\frq:M\pol N$
 is a measure $\frq$ on $M\times N$
such that the pushforward%
\footnote{Let $(M,\mu)$ be a measure space, let $S$ be a set, $\psi:M\to S$ a map. We say that a subset $C\subset S$
is measurable if $\psi^{-1}(C)$ is measurable, the pushforward $\sigma$ of the measure $\mu$
is defined by $\sigma(C)=\mu (\psi^{-1}(C))$.}
of $\frq$ under the projection $M\times N\to M$ coincides
 with $\mu$ and the pushforward under the projection $M\times N\to N$ coincides
 with $\nu$. 
 
 \sm 
 
 {\sc Example.} Let $g\in \Ams(M)$. Consider the map 
 $M\to M\times M$ given by $m\mapsto (m, mg)$. Then the pushforward 
 of $M$ under this map is a measure preserving polymorphism.
 In particular, for the identical transformation $g=1$ we get an {\it identical polymorphism.}
 \hfill $\boxtimes$
 
 \sm
 
 A measure preserving polymorphism $\frq:M\pol N$  can be regarded as a 'spreading
 map' sending points to measures. Namely, consider a partition
  $$M\times N=\coprod\limits_{m\in M} m\times N$$ 
  and conditional measures%
  \footnote{See \cite{Roh}, \cite{Bog}, Chapter 10.} $\kappa_m(n)$ on $N$ defined for almost
  all $m$ from the condition
  $$
  \frq(S)=\int_M \kappa_m\bigl((m\times N)\cap S\bigr)\, d\mu_m(n)
  \quad \text{for any measurable $S\subset M\times N$.}
  $$
Then {\it we can regard $\frq$ as a 'spreading map' sending points $m\in M$ to measures
$\kappa_m(n)$ on $N$}.

A {\it product} of measure preserving polymorphisms corresponds to a double spreading.
Namely, let $\frq$ be a polymorphism $M\pol N$ and $\frp$
be a polymorphism $N\pol K$. Let $\kappa_m(n)$ be the system of conditional measures corresponding to $\frq$, and $\pi_n(k)$
be the system of conditional measures corresponding to $\frp$. 
Then the system $\rho_m(k)$ of conditional  measures
corresponding to the product $\frr=\frp\circledast \frq$
is defined by
$$
\rho_m(k)=\int_N \pi_n(k) \,   d\kappa_m(n).
$$

The map $g\mapsto T(g)$, see (\ref{eq:T-pres}), can be extended from $\Ams(M)$ to polymorphisms $\frq:M\pol N$. Namely,
 the  {\it operator}
$T(\frq):L^2(N)\to L^2(M)$ is given by
$$
T(\frq) f(m)=\int_N f(n) \,d\mu_m(n).
$$
Then for $\frq:M\pol N$, $\frp:N\pol K$ we have
$$T(\frp\circledast \frq)=T(\frp)\, T(\frq).$$

 \sm 

We say that a sequence $\frq^{(j)}:M\pol N$ {\it converges} 
to $\frq$ if for any measurable $A\subset M$, $B\subset N$
the sequence $\frq^{(j)}(A\times B)$ converges to
$\frq(A\times B)$. This convergence is equivalent to the weak operator convergence 
$T(\frq^{(j)})\to T(\frq)$. 

The product of polymorphisms is separately continuous;
the group $\Ams(M)$ is dense in the semigroup of polymorphisms  $M\pol M$.

\sm 

We also have the {\it involution} $\frq\mapsto \frq^\star$ on polymorphisms, it sends a measure $\frq$ on $M\times N$ to the
same measure on $N\times M$. We have $T(\frq^\star)=T(\frq)^*$.

\sm

%For $g\in \Ams(M)$ we consider a map $M\mapsto M\times M$
%defined by $m\mapsto (m,g(m))$. The image $\frQ(g)$ of the measure $\mu$ under this map is a bistochastic kernel. Then $T(\frQ(g))=T(g)$ and 

 \sm 
 
 {\bf \punct Quotient spaces and the corresponding measure preserving polymorphisms.%
 \label{ss:quotients}}
 Let $M$ be a Lebesgue measure space, let $M=\sqcup_{h\in H} S_h$ be a measurable partition,
 let $\sigma_h(m)$ be the conditional measures on sets $S_h$
 (see \cite{Roh}).  Consider the quotient space $H$ with the induced measure $\eta(h)$.
 Then we have the following polymorphisms (see \cite{Ner-boundary}, Subsect. 3.10): 
 
 \sm 
 
 a) We have the projection  $\pi:M\to H$ sending $m\in M$ to the element $S_h$ containing $m$.
 We consider the map $M\to M\times H$ given by $m\mapsto (m,\pi(m))$, we denote by $\frx$
 the image of $\mu$ under this map. Then $\frx$ is a measure-preserving polymorphism  $M\pol H$.
 The operator
 $$
 T(\frx)f(m)=f(\pi(m))
 $$
  is a canonical isometric embedding $L^2(H)\to L^2(M)$. The image consists of $L^2$-functions that are constant on sets $S_h$. 
 
 \sm 
 
 b) We also have a polymorphism $\frx^\star$, the corresponding operator $T(\frx^*):L^2(M)\to L^2(H)$
 is the operator of {\it conditional expectation}
 $$
 T(\frx^*)f(h)= \cE f(h):=\int_{S_h} f(m)\,d\sigma_h(m).
 $$

 \sm
 
 c) Notice that $\frx  \circledast \frx^\star$ is the identical polymorphism $H\to H$.
 
 \sm
 
 d) For the polymorphism  $\frx^\star \circledast \frx:L^2(M)\to L^2(M)$
 the corresponding operator is the orthogonal projector   to the subspace $L^2(H)$.
 
 \sm
 
 {\bf \punct General polymorhisms.%
 \label{ss:poli}}
 Let $M$ be a Lebesgue space with a continuous  measure.
Recall that a map $g:M\to M$ {\it leaves a measure $\mu$ quasiinvariant}%
\footnote{This is also equivalent to the condition: $g$ is bijective a.s. and both $g$ and $g^{-1}$ sent sets of zero measure to sets of zeto measure.}
if it  is bijective almost sure and there is a function $g'(m)$ 
(the {\it Radon--Nikodym derivative}, see, e.g., \cite{Bog}, 9.12) such that for any measurable
$A\subset M$ we have
$$
\mu(Ag)=\int_A g'(m)\,d\mu(m),
$$
the $g'(m)$ is a natural extension of Jacobians in classical analysis.
Denote by $\Gms(M)$ the group of transformations of $M$ leaving the measure $\mu$ quasiinvariant.

\sm 

Denote by $\Pi\subset \C$ the strip 
 \begin{equation}
 0\le \Im \lambda\le 1.
 \label{eq:Pi}
 \end{equation}
 For $r+is\in \Pi$
  we consider the following transformations
of the space of functions on $M$:
\begin{equation}
T_{ r+is}(g)\, f(m)=f(mg)\,g'(m)^{r+is},
\label{eq:p+is}
\end{equation}
 Then $T_{ r+is}(g)$ 
is an isometric operator $L^{1/r}(M) \to L^{1/r}(M)$. 

\sm

{\sc Remark.}
The main topic of our
interest is the case $r=1/2$, when we get a unitary operators $L^2\to L^2$. However, operators in spaces $L^p$ are  used in our argumentation in Theorem \ref{th:1} below.
\hfill $\boxtimes$

\sm

Denote by $\R^\circ$ the multiplicative group of positive numbers, denote by $t$ the coordinate on it.
Let $(M,\mu)$, $(N,\nu)$ be probabilistic Lebesgue spaces. We say that a {\it polymorphism} $\frp: M\pol  N$ is a measure on $M\times N\times \R^\circ$
such that 

\sm

1. The pushforward of $\frp$ under the projection 
$M\times N\times \R^\circ\to M$ is the measure $\mu$.

\sm 

2. The pushforward of the measure $t\cdot \frp$ to $N$ coincides with the measure
$\nu$.

\sm 

Denote by $\Pol(M,N)$ the space of all polymorphisms $M\pol N$.
%By $\Pol(M)$ we denote the space of all polymorphisms from $M$ to itself.
See \cite{Ner-bist}, \cite{Ner-boundary}.

\sm 

{\sc Example.} Measure preserving polymorphisms  from Subsect. \ref{ss:bist} 
are polymorphisms, they are supported by the set
$M\times N\times \{1\}\subset M\times N\times \R^\circ$.
\hfill $\boxtimes$

\sm

{\sc Example.} Let $q\in \Gms(M,\mu)$. Consider the map 
$M\to M\times M\times \R^\circ$ defined by
$$
m\mapsto \bigl(m,q(m), q'(m)\bigr).
$$
Then the pushforward of the measure $\mu$ under this map is
a polymorphism $M\pol M$.
\hfill $\boxtimes$

\sm 

{\sc Example.}
Denote by $\cM^\triangledown(\R^\circ)$ the set of all positive finite Borel measures $\sigma$ on $\R^\circ$
such that $t\cdot \sigma$ also is finite. Clearly, any measurable function
$(m,n)\mapsto s_{m,n}$
from $M\times N$ to $\cM^\triangledown$ determines a certain measure $\frs$
on $M\times N\times \R^\circ$. Namely, for measurable sets
$A\subset M$, $B\subset N$, $C\subset \R^\circ$
we set
$$
\frs(A\times B \times C)=\int_{A\times B} s_{m,n}(C)\,d\mu(m)\,d\nu(n).
$$
We also must assume that
$$
\forall m:\,
\int_{N} s_{m,n}(\R^\circ)\,d\nu(n)=1,
\qquad
\forall n:\,
\int_M \int_{\R^\circ} t\, ds_{m,n}(t)\, d\mu(m)=1,
$$
under these conditions $\frs $ is contained in $\Pol(M,N)$.
We say that a polymorphism $\frs$ obtained in such  way is {\it continuous}. 
Notice that polymorphisms from the previous example are not continuous.
\hfill $\boxtimes$

\sm 

{\bf \punct Convergence of polymorphisms.}
Let $\frp_j\in\Pol(M,N)$. Consider measurable subsets $A\subset M$, $B\subset N$. Restrict $\frp$ to $A\times B\times \R^\circ$ and take its
pushforward under the map $A\times B\times \R^\circ\to\R^\circ$.
Denote the resulting measure on $\R^\circ$ by $\frp[A\times B]$. 
We say that a sequence  $\frp_j\in \Pol(M,N)$ converges to $\frp$, 
if for any measurable subsets $A\subset M$, $B\subset N$ we have weak
convergences of measures
$$
\frp_j[A\times B]\to \frp[A\times B],\qquad  t\cdot \frp_j[A\times B]\to t\cdot\frp[A\times B].
$$

It is easy to show (see \cite{Ner-boundary}, Theorem 5.3) that {\it the group $\Gms(M)$ is dense in $\Pol(M,M)$}.
Continuous polymorphisms $M\pol N$  are dense in $\Pol(M,N)$.

\sm 

{\bf \punct Products of polymorphisms.}
There is an obvious way to multiply continuous polymorphisms. Let
$\fru\in \Pol(M,N)$, $\frv\in \Pol(N,K)$ be determined by functions
$(m,n)\mapsto u_{m,n}$, $(n,k)\mapsto v_{n,k}$.
Then their product 
$$\frw=\frv \circledast\fru \in \Pol(M,K)$$
is determined by the function
$M\times K\to \cM^\triangledown$ given by
$$
w_{m,k}=\int_N u_{m,n}*v_{n,k}\, d\nu(n),
$$
where $*$ denotes the convolution on $\R^\circ$.

{\it The multiplication $\circledast$ extends to a separately continuous map}
	$$
	\Pol(M,N)\times \Pol(N,K)\to \Pol(M,K) ,
	$$
and moreover this operation is associative (see \cite{Ner-boundary}, Theorem 5.5, Theorem 5.9). So we get a category
$\Pol$ whose objects are Lebesgue probabilistic measure spaces and morphisms are polymorphisms. 

For different ways to define the multiplication of polymorphisms, see \cite{Ner-boundary}, however it seems that definitions by continuity
(by continuous extension from continuous polymorphisms or from the group $\Gms(M)$) are
more convenient.

\sm 

{\bf \punct Polymorphisms as spreading maps.%
\label{ss:point-image}}
 Informally, {\it polymorphisms are 'maps' $M\to N$, which spread points of $M$ along $N$ and the Radon--Nikodym derivative also is spread.}  Namely, let $\frq\in \Pol(M,N)$.
 Consider a partition 
 $$
 M\times N\times \R^\circ=\coprod\limits_{m\in M} m\times N\times \R^\circ
 $$
 and conditional measures $\kappa_m(n,t)$
 on fibers of the partition. We say that   {\it measures $\kappa_m(n,t)$ are spreaded images of points $m$}.
 Below in Section 3
 we describe polymorphisms of $\frS^z$ in terms of such conditional measures.

\sm 

{\bf \punct Involution.}
Define  an involution in the category $\Pol$.
Consider the map $(m,n,t)\mapsto (n,m,t^{-1})$.
For $\frp\in \Pol(M,N)$ we consider its pushforward $\frp'$ under this map and  define $\frp^\star\in \Pol(N,M)$ as the measure $t\cdot \frp'$. Then 
$$
(\frq\circledast \frp)^\star=\frp^\star \circledast \frq^\star.
$$

\sm 

{\bf \punct Mellin--Markov transforms of polymorphisms.}
For each $r+is\in \Pi$
the map $g\mapsto T_{r+is}(g)$ extends  to the category of polymorphisms.

Let $\frp\in \Pol(M,N)$. The {\it Mellin--Markov transform of $\frq:M\pol N$} 
is a function 
defined in the strip $r+is\in \Pi$ taking values
in operators
$$
T_{r+is}(\frq):L^{1/r}(N)\to L^{1/r}(M)
$$
and determined
 from the equality
\begin{equation}
\int\limits_M \phi(m)\,T_{r+is}(\frq) \psi(m)\,d\mu(m)=
\iint\limits_{M\times N \times \R^\circ}
\phi(m) \psi(n) t^{r+is} \,d\frq(m,n,t)
\end{equation}
for any $\psi\in L^{1/r}(N)$, $\phi\in L^{1/(1-r)(M)}$.

The operators $T_{r+is}(\frq):L^{1/r}(N)\to L^{1/r}(M)$ are bounded and
 for any $\frp\in\Pol(M,N)$, $\frq\in\Pol(M,K)$ we have (see \cite{Ner-boundary}, Theorem 6.14)
$$
T_{r+is}(\frq)\, T_{r+is}(\frp)=T_{r+is}(\frq\circledast \frp).
$$

It is more transparent to determine
 $\frq$  in terms of conditional measures $\kappa_n(m,t)$ as in Subsect. \ref{ss:point-image}. Then
$$
T_{r+is}(\frq) f(m)=\int_{N\times \R^\circ}
 f(n)\,t^{r+is}\,d\kappa_m(n,t).
$$

%First, let $\fru\in\Pol_\cont (M,N)$ be defined by a function $u:M\times N\to \cM^\triangledown$.
%Then an operator
%$$
%T_{r+is} f(m)=\int_N \int_{\R^\circ}  t^{r+is} f(n) \, du_{n,m}(t)\, d\nu(n) 
%$$
%sends functions on $N$ to functions on $M$. This operator is bounded as an operator
%$L^{\frac 1p}(M)\to L^{\frac 1p}(N)$ and, moreover, is {\it contractive} has norm $\le 1$. More generally, for an arbitrary 
%$\fru\in \Pol(M,N)$

%For $\frp\in\Pol(M,N)$ we define a
%Hermitian form on $L^{\frac 1{1-p}}(M)\times L^{1/p}(N)$ by
%$$
%\{\phi,\psi\}=\int_{M\times N\times \R^\circ}  t^{r+is} \phi(m)\,\ov{\psi(n)}\,d\fru(m,n,t).
%$$
%It can be shown that

%$$
%\bigl|\{\phi,\psi\}\bigr|\le \|\phi\|_{\frac 1{1-p}(M)} \cdot \|\psi\|_{L^p(N)}.
%$$
%By duality this define the operator
%$$
%T_{p+is}: L^{\frac 1p} (N)\to L^{\frac1p}(M)
%$$
%with norm $\le 1$. Moreover, 
%$$
%T_{p+is}(\frv\circledast \fru)=
%T_{p+is}(\frv)\,T_{p+is}(\fru)
%$$

%Dzyaguguk12345
\sm

%{\sc Example. A special case: measure preserving polymorphisms.}
%Denote by $\Pol_1(M,N)$ the space of polymorphisms $M\pol N$ supported by $M\times N\times 1\subset M\times N\times\R^\circ$. So we have measures $\fru$ on $M\times N$ whose pushforwards
%to $M$ and $N$ are $\mu$ and $\nu$ respectively. Product of such polymorphisms has the same
%form, in this case there are more transparent definitions of product, see \cite{Kre}, \cite{Ver}, \cite{Ner-book}.
%\hfill $\square$

{\bf \punct Closures of actions of groups on measure spaces.}
For a group $G$ acting on a measure space $M$ we get a question about the closure
of $G$ in the semigroup $\Pol(M,M)$ of  polymorphisms. This question for $G=\Z$ acting 
by measure preserving transformations was a subject of numerous works, see e.g., \cite{King}, \cite{JPR}, \cite{Sol}, \cite{KR},
a related topic is
 intertwiners (joinings)  of such actions   in class of measure preserving polymorphisms.

For infinite-dimensional groups
the first problem of this kind was solved by Nelson \cite{Nel}, who examined  the action
of the infinite-dimensional orthogonal group on the infinite-dimensional space with a Gaussian measure (the resulting semigroup is the semigroup of all contractive operators).
For the natural group of symmetries of Poisson measure and for the action of the 
restricted group $\GL$ on a space with a Gaussian measure the problems were discussed in \cite{Ner-match}, \cite{Ner-gauss}. They are unexpectedly non-trivial and lead to formulas, which are at least unusual. One simple case (measures on the space of infinite Hermitian matrices invariant with respect to unitary groups) was examined in \cite{Ner-Whi}.

In this paper we consider the action of the bisymmetric group $\SS$ on spaces $\frS^z$ of virtual permutations.

\sm 

{\bf \punct The purpose of the paper.}
First, we show that the categories $\cS$ and $\ov\cS$
  of chips act on the space of virtual permutations by polymorphisms.
The category of chips is a representative of a wide zoo of train constructions,
 in Section 2 we prove a general statement in the following spirit: certain actions of infinite-dimensional group  on  measure spaces
 generate actions of their trains by 
polymorphisms.

In Section 3 we write a formula for the action of the semigroup $\cS(0,0)$ on the space of $\frT^z$ of restaurants. In fact, for a given collection $\bigl\{\bigcirc[k_j]\}$ of cycles we write spreaded images of points as  sums over a certain sets of dessins d'enfant (or equivalently of checker triangulated surfaces), summands are
 multiple convolutions of Dirichlet distributions (Theorem \ref{th:2}).
 
 In Section 4 we get some further statements in this spirit.  In Theorem \ref{th:3} we get formulas for the action of the center
 $\cS(0,0)$
 of the semigroup $\frS(\infty,\infty)$ on $\frS^z$. 
 
 In Theorem \ref{th:2} we write a similar formula for arbitrary chips, whose left sides are trivial (see Fig.~\ref{fi:left-trivial}).
 Formula involves several types of parameters of combinatorial nature and is longer than for the action of $\cS(0,0)$
 This allows to write explicit formulas in the same spirit for generators of the category $\cS$ (we omit this).
 However, the author did not succeed to obtain a general formula for arbitrary chips.

\section{Trains of $(G,K)$-pairs and polymorphisms.}

\COUNTERS

Train constructions (multiplication of double cosets) and the multiplicativity theorem
\ref{th:0} are relatively usual phenomena for infinite-dimensional groups, including infinite symmetric groups \cite{Olsh-symm}, \cite{Ner-imrn}, \cite{Ner-umn}, 
classical groups over reals \cite{Olsh-GB}, \cite{Ner-book},
\cite{Ner-colligations} (and also over finite fields
\cite{Olsh-semi}, \cite{Ner-finite} and $p$-adic fields \cite{Ner-p}, \cite{Ner-Bruhat})
 groups of transformations of measure spaces
 \cite{Ner-bist}, \cite{Ner-book}  and some exotic cases.
 
 For infinite-dimensional group   a lot of 
 quasiinvariant actions (i.e., embeddings to the group $\Gms(\cdot)$) is known. 
 In this section, we wish to show that under some conditions  such actions automatically generate actions of corresponding
 train categories by polymorphisms. We prove Theorem \ref{th:1},
 which apparently is sufficient
to justify this claim for known zoo of train constructions. Also, the proof is  simple, short,
and can be easily repeated independently in any explicit case.

We do not pretend to any wider purpose.

\sm 

{\bf\punct Operators corresponding to double cosets.}
Let $G$ be a separable topological group, let $K$ be a closed subgroup. Let $\rho$ be a unitary representation 
of $G$ in a separable Hilbert space $H$. Denote by $H(K)$ the subspace
consisting of vectors fixed by elements of $K$, by $P(K)$ the orthogonal
projection to $H(K)$. For two subgroups $K$, $L$ and $g\in G$ consider
the operator
$$
\wt\rho_{(K,L)}(g):=P(K)\rho(g)\bigr|_{H(L)}:\,H(L)\to H(K).
$$
It is easy to verify that the operator-valued function $\wt\rho_{(K,L)}(g)$ is constant on double
cosets $KgL$, so it is a function on $K\backslash G/L$.
 We define a more 
coarse equivalence relation on $G$,
\begin{equation}
g\sim g' \quad \text{if} \quad \wt\rho_{(K,L)}(g)=\wt\rho_{(K,L)}(g')\quad
\text{for all $\rho$}.
\label{eq:coarse}
\end{equation}
Denote an equivalence class containing $g$ by $[K g L]$,
denote the space of equivalence classes by $[K\backslash G/L]$.

\sm 

{\bf\punct $(G,K)$-pairs.}
Now let $G$ be a separable topological group, let $K$ be a subgroup.
Let $\{K_\alpha\}_{\alpha \in \cA}$ be
a family of subgroups in $K$. Let the set $\cA$ be partially ordered,
and for any $\alpha$, $\alpha'\in A$ there is $\nu\in \cA$ such that 
$\nu \succcurlyeq \alpha$, $\nu \succcurlyeq \alpha'$. Assume also
that $\alpha \succcurlyeq \alpha'$ implies  $K_{\alpha'}\supset K_{\alpha}$.

\sm 

A {\it $(G,K)$-pair} is a group and a family of subgroups $\{K_\alpha\}_{\alpha\in\cA}$ as above satisfying the following conditions: 

\sm 

A. {\it For each $\alpha$, $\beta$, $\gamma\in \cA$, for each
$\frg\in [K_\alpha\bs G/K_\beta]$, $\frh\in [K_\beta\bs G/K_\gamma]$
there is $\frp\in [K_\alpha\bs G/K_\gamma]$ such that for each unitary representation $\rho$ of $G$ in a Hilbert space $H$ we have
\begin{equation}
\wt\rho_{(K_\alpha,K_\beta)}(\frg)\,\,
 \wt\rho_{(K_\beta, K_\gamma)}(\frh)=
\wt\rho_{(K_\beta,K_\gamma)}(\frp).
\end{equation} }

B. {\it For each unitary representation, the subspace $\cup_{\alpha\in\cA} H(K_\alpha)$ is dense in $H$.}

\sm 

So we get an operation
$$
[K_\alpha\bs G/K_\beta]\times [K_\beta\bs G/K_\gamma]\,\to \,[K_\alpha\bs G/K_\gamma].
$$
Denote it by $\frp=\frg\circ \frh$ or $\frp=g\circ h$.
Since the product of linear operators is associative, this operation
is also associative. So we get a category -- {\it the train $\cT=\cT(G,K)$ of $(G,K)$}.
Its objects are enumerated by the set $\cA$,  morphisms $\beta\to\alpha$
are reduced double cosets $[K_\alpha\bs G/K_\beta]$.

 This category is equipped with 
the {\it involution} $\frg\mapsto\frg^\ast$, which is the map 
$$[K_\alpha\bs G/K_\beta]\to [K_\beta\bs G/K_\alpha]$$
induced
by the inversion in $G$, $g\mapsto g^{-1}$.

To simplify notations, denote
$$
H[\alpha]:=H(K_\alpha),\quad P[\alpha]:=P(K_\alpha),
\quad 
\wt \rho_{\alpha,\beta}:=\wt\rho_{(K_\alpha,K_\beta)}.
$$

\sm 

{\bf \punct Actions of $(G,K)$-pairs by polymorphisms.}
Consider a $(G,K)$-pair. Let $G$ act on a Lebesgue probabilistic measure space $(M,\mu)$ by transformations leaving the measure quasiinvariant, let $K$ act by measure preserving transformations.
In such a case we say that we have an {\it action of a $(G,K)$-pair
	on a measure space.}

For each $K_\alpha$  consider the $\sigma$-algebra 
$\Sigma[\alpha]$ of all $K_\alpha$-invariant sets. Namely, a set
 $B\subset M$ is contained in $\Sigma[\alpha]$ if it is measurable and for any $g\in G$ the symmetric sum
$B\triangle Bg$ has measure 0; equivalently, the indicator function of
$B$ is a $K_\alpha$-fixed element of $L^2(M,\mu)$.
This sigma-algebra $\Sigma[\alpha]$ determines a measurable partition%
\footnote{Recall that a partition $M=\sqcup_{h\in H} S_h$ of $M$ is {\it measurable} if there is a countable family 
of measurable subsets $A_j$, which are unions  of elements of partitions, such that for each distinct elements  $S_{h}$, $S_{h'}$ there is $A_j$
such that $S_h\subset A$,
$S_{h'}\subset M\setminus A$ or vise versa. 
\newline 
Consider 
the set $\cI[\alpha]\subset L^2(M)$ of all indicator functions $I_C$ of elements of $\Sigma[\alpha]$. Choose a dense countable 
subset $\{I_{C_j} \}$ in $\cI[\alpha]$. The sigma-algebra generated by $C_j$ and sets of zero measure coincides with $\Sigma[\alpha]$.
We say that $m$, $m'\in M$
are equivalent if for each  $C_j$ we have $m$, $m'\in C_j$
or $m$, $m'\notin C_j$. By the definition, we get a measurable partition. By the Rohlin theorem \cite{Roh} the quotient is a Lebesgue measure space.}
of $M$, the quotient space is a Lebesgue measure space, denote them  $\bigl(M[\alpha],\mu[\alpha]\bigr)$. We have an identification
of $L^2\bigl(M[\alpha]\bigr)$ with the space of $K_\alpha$-fixed vectors
in $L^2(M)$.

By Subset. \ref{ss:quotients},  %(see \cite{Ner-boundary}, Subsect. 3.10)
we get a canonical measure preserving polymorphism $\frx[\alpha]: M\pol M[\alpha]$. 
%Namely, we have a map
%sending each point $m\in M$ to the element $\ov m$ of the partition
% containing $m$. Therefore, we have a map $M\to M\times M[\alpha]\times \R^\circ$
% given by $m\mapsto (m,\ov m,1)$. We define $\frx[\alpha]$ as the pushforward of $\mu$ under this map.
 
 Operators  $T_{r+is}\bigl(\frx[\alpha] \bigr): L^{1/r}\bigl(M[\alpha]\bigr)\to L^{1/r}(M)$
 are  operators
 $$
 \cE_{1/r}'[\alpha]\, \phi(m)=\phi(\ov m),\quad\text{where $\ov m$ is  element of the parition containing $m$,}
 $$
 these operators do not depend on $s$ and in fact do not depend on $r$.
 
 An operator
 $T_{r+is}\bigl(\frx[\alpha] \bigr)$ 
 is the operator of conditional expectation
 %(see, e.g.,\cite{Bog},) 
 $\cE_{1/r}[\alpha]:L^{1/r}(M)\to L^{1/r}\bigl(M[\alpha]\bigr)$.
For $r=1/2$,  $T_{1/2+is}\bigl(\frx[\alpha]^\ast \bigr)$  is the orthogonal
 projection to $L^2\bigl(M[\alpha]\bigr)$.

 %An operator $T_{1/2+is}\bigl(\frl[\alpha] \bigr)$
 %does not depend on $s$.
 %More generally, 

 %embedding $L^2\bigl( M[\alpha]\bigr )\to L^2(M)$. The adjoint operator

 \begin{theorem}
 	\label{th:1}
 	Let a $(G,K)$-pair act on a measure space $(M,\mu)$.
 	Assign for each $\alpha\in \cA$ the measure space $M[\alpha]$.
 	For each $\alpha$, $\beta\in\cA$ and $g\in G$
 	 define the polymorphism
 	\begin{equation}
 	\Xi_{\alpha,\beta}(g):= \frx[\alpha]^\star \,g\, \frx[\beta] \,\in \Pol\bigl(M[\alpha], M[\beta]\bigr).
 	\label{eq:xi}
 	\end{equation}
 	Then the maps $g\mapsto \Xi_{\alpha,\beta}(g)$ determine
 	a functor from the train $\cT(G,K)$ to the category
 	of polymorphisms.
 \end{theorem}

{\sc Proof.} 
We must verify   the identity
\begin{equation}
\Xi_{\alpha,\beta}(\bigl[K_\alpha g K_\beta]\bigr)\circledast \Xi_{\beta,\gamma}\bigl ([K_\beta h K_\gamma]\bigr)
=\Xi_{\alpha,\gamma}(\bigl ([K_\alpha g K_\beta]\circ [K_\beta h K_\gamma] )
\label{eq:prove}
\end{equation}
By \cite{Ner-boundary}, Theorem 6.12, two polymorphisms are equal 
if and only if their Mellin--Markov transforms are equal. A Mellin--Markov transform
is a holomorphic operator-valued function in the strip $\Pi$, see \cite{Ner-boundary}, Lemma 6.10, so it is uniquely determined
by its values on the line $1/2+is$.
Applying functors $T_{1/2+is}$ to both sides of conjectural equality (\ref{eq:prove})
we come to the equivalent conjectural equality
\begin{equation*}
\Bigl(\cE[\alpha] T_{1/2+is}(g)\cE'[\beta]\Bigr)\, \Bigl(\cE[\beta]\, T_{1/2+is}(h)\cE'[\gamma]\Bigr)
=\cE[\alpha]  T_{1/2+is}(g\circ h)\cE'[\gamma],
\end{equation*}
these operators act $L^2(M[\gamma])\to L^2(M[\alpha])$.
Since $\cE'[\beta]\cE[\beta]=P[\beta]$, we can write the 
left hand side as
$$
P[\alpha]\,T_{1/2+is}(g)\, P[\beta]\, T_{1/2+is}(h)\Bigr|_{L^2(M[\gamma])}.
$$
By the definition of trains, this equals to
$$
P[\alpha]\,T_{1/2+is}(g\circ h)\Bigr|_{L^2(M[\gamma])},
$$
i.e., to the right-hand side.
\hfill $\square$

\sm

{\bf \punct Application to the bysymmetric group.} 
1) We  take $G=\SS$, $K[\alpha]=\K_\alpha$, $M=\frS^z$,
and get an action of the category $\cS$.

\sm

2) In our case $[\K_\alpha\backslash \SS /\K_\beta]=\K_\alpha\backslash \SS /\K_\beta$, this is  clear from
explicit  constructions of representations of $\SS$, see \cite{Olsh-symm}. However, in  considerations above
sets $[K_\alpha\backslash G/K_\beta]$ 
are  used
only for the establishing of the assiciativity, in our case the associativity is obvious.
So we can simply repeat the proof of Theorem \ref{th:1} without a reference to 
\cite{Olsh-symm}.

\sm

3) If we want to consider the extended category $\ov\cS$, see Subsect. \ref{ss:infty-infty}, then we
take a countable subset $\Omega\subset \N$ such that $\N\setminus \Omega$ also is countable.
Elements of  $\cA$ are subsets $C\subset \N$ such that $\Omega\setminus C$ are finite.
A group $K(C)$
is the subgroup in $\K\simeq \ov S_\infty$ fixing all elements of $C$.

%%%%%%%%%%%%%%%%%%%%%%%%%%%%%%%%%%%%%%%%%%%%%%%%%%%%%%%%%%

%%%%%%%%%%%%%%%%%%%%%%%%%%%%%%%%%%%%%%%%%%%%%%%%%%%%%%%%%%

%%%%%%%%%%%%%%%%%%%%%%%%%%%%%%%%%%%%%%%%%%%%%%%%%%%%%%%%%%

%%%%%%%%%%%%%%%%%%%%%%%%%%%%%%%%%%%%%%%%%%%%%%%%%%%%%%%%%%

%%%%%%%%%%%%%%%%%%%%%%%%%%%%%%%%%%%%%%%%%%%%%%%%%%%%%%%%%%

\section{Action of the semigroup $\cS(0,0)$ on the space $\frT^z$ of restaurants}

\COUNTERS

In the construction of the previous section we set:
$$M=\frS^z,\qquad G=\SS,\qquad K=\K, \qquad K_n=\K_n$$

\sm

{\bf \punct Spaces of half-empty restaurants.} Denote by $\frS_n^z$ the measure space, whose points  are restaurants  $\{U_\omega\}\in\frT^z$ equipped with guests
$\ov 1$, \dots, $\ov n$ chosen uniformly,  denote such points by $\{U_\omega\},\{\ov j\}_{j\in I_n}$. 
Denote the measure on $\frS_n^z$ by $\ov\mu_n^z$. In particular $\frS_0^z=\frT^z$.
%We call lengths of intervals between neighbor guests by {\it interludes}. %[vstavka, birka]
%To uniformize terminology, we also call by interludes lengths of empty tables. 

We have an obvious forgetting map $\frS^z\to\frS^z_n$, the pushforward of the measure $\mu^z$ is $\ov \mu^z_n$.

\begin{proposition}
$$
\frS^z[n]= \frS_n^z.
$$ 
\end{proposition}

{\sc Proof.}
Since $\frS_n^z$ is the quotient space of $\frS^z$, the space $L^2(\frS^z_n)$ is a canonically defined subspace
in $L^2(\frS^z)$ and elements of this subspace are $\K$-fixed. 
 We must show that the space of $\K_n$-fixed vectors in $L^2(\frS^z)$ is precisely the $L^2(\frS^z_n)$.
Notice that $L^2(\frS^z)$ is a direct integral of  Hilbert spaces,
 the base of the integral is the measure space $\frS^z_n$
 and a fiber over a point $\{U_\omega\},\{\ov j\}_{j\in I_n}$  is $L^2$ on a  product of countable number of copies of
 $\sqcup U_\omega$ enumerated by $n+1$, $n+2$, \dots.
  The group $\K_n$ acts trivially on the base $\frS^z_n$. In each fiber 
 it acts by permutation of factors $\sqcup U_\omega$. So a $\K_n$-fixed vector in $L^2(\frS^z)$
is an integral of $\K_n$-fixed vectors in fibers $L^2\bigl((\sqcup U_\omega)^\infty\bigr)$. By the zero-one law such a vector is a constant function.
\hfill $\square$

\sm

{\sc Remark.} So the polymorphism $\frx[n]^\star:\frS^z\to \frS^z_n$ is the removing
of guests $\ov {n+1}$, $\ov {n+2}$, \dots. The adjoint polymorphism  $\frx[n]:\frS^z_n\to \frS^z$
is a random arrangement of new guests $\ov {n+1}$, $\ov {n+2}$, \dots
(these polymorphisms are measure-preserving).

\sm 

\begin{figure}
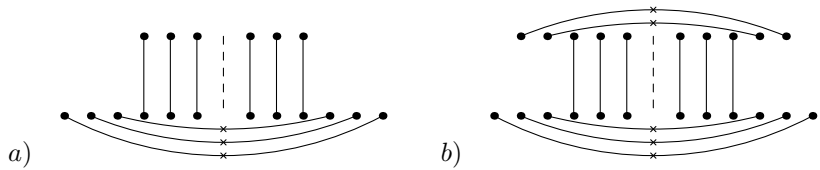

	$$a)\quad \epsfbox{chips-add.1}
	\qquad b)\quad 
	\epsfbox{chips-add.2}   
	$$
	\caption{a) A chip $\lambda^6_3$.
	\newline 
b) To Subsect. \ref{ss:embeddings}.}
\label{fi:lambda}
	\end{figure}

{\bf \punct Embeddings and projections.%
\label{ss:embeddings}} Let $m\ge n$.
Define the following canonical morphisms $\lambda^m_n\in \cS(n,m)$, see Fig.~\ref{fi:lambda}.a. We take arcs

\sm 

1) $1_r^+\downarrow 1_r^-[0]$, \dots, $n_r^+\downarrow n_r^-[0]$ in the right hand side;

\sm 

2) $1_l^+\downarrow 1_l^-[0]$, \dots, $n_l^+\downarrow n_l^-[0]$ in the left hand side;

\sm 

3) $(n+1)_l^+\smallsmile (n+1)_r^+[1/2]$,\, \dots,\, $m_l^+ \smallsmile m_r^+[1/2]$. 

\sm 

The following statement immediately follows from the definition:

\begin{proposition}
	\label{pr:lambda}
	{\rm a}  The polymorphism $\Xi_{m,n}\bigl((\lambda^n_m)^*\bigr): \frS^z_n\to \frS^z_m$
	is a forgetting of guests $\ov{m+1}$, \dots, $\ov n$.
	
	\sm 
	
	{\rm b)}  The polymorphism $\Xi(\lambda^m_n): \frS^z_m\to \frS^z_n$ is a uniform random arrangement
	of guests $\ov{m+1}$, \dots, $\ov n$ on a given restaurant $\sqcup U_\omega$.
\end{proposition}

More formally, 
let
$\bigl(\sqcup U_\omega, \{\ov i\}_{i\le n}\bigr)\in \frS^z_n$.
Then its spreaded image  under $\frx(\lambda_m^n)$
is the measure on $\R^\circ\times \frS^z_m$ supported by the set consisting of points
$$1\times\bigl(\sqcup U_\omega, \{\ov i\}_{i\le n}, 
\{\ov j \}_{n<j\le m}\bigr) ,$$
 where the guests $\{\ov i\}_{i\le n}$ are the same, and 
guests $\ov{n+1}$, \dots, $\ov m$ are uniformly distributed on
$\sqcup U_\omega$.

\sm 

{\sc Remark.} Proposition \ref{pr:lambda} also
provides us a description of polymorphisms corresponding to diagrams drawn on Fig.~\ref{fi:lambda}. We forget part of guests
and arrange new guests randomly.
\hfill $\boxtimes$

\sm 

%Consider slightly more general example.
%Next, fix $k$ and $n$, $m\ge k$. Consider a chip (see Fig.
%\ref{fi:lambda}.b) consisting
%of arcs

%\sm

%a) $j_l^+\downarrow j_l^-$, $j_r^+\downarrow j_r^-$ for $j\le k$;

%\sm 

%b) $i_l^+\smallsmile i_r^+$ for $i=k+1$, \dots, $n$;

%\sm 

%c) $i_l^-\smallfrown i_r^-$ for $i=k+1$, \dots, $m$.

%\sm 

%Then the corresponding polymorphism.

%Consider the map 
%$$\frS_n^z\times \frS_m^z\to \frS_k^z\times \frS_k^z$$

\sm 

Our next purpose is to describe
 the action of the semigroup $\cS(0,0)$ on the space $\frT^z=\frS_0^z$.
We need to define some additional objects, namely checker surfaces (they also arise in representation theory of infinite symmetric group
for other reasons, see \cite{Ner-umn},\cite{Ner-imrn}).

\sm

{\bf\punct Checker surfaces.} See \cite{Ner-imrn}, \cite{Ner-umn}.  Consider an oriented  compact two-dimensional closed surface (generally speaking it is disconnected).
\begin{figure}
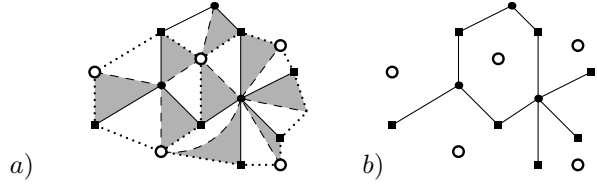

 $$
a)\qquad \epsfbox{ribbon.1}  \qquad b) \epsfbox{ribbon.2} 
$$
\caption{
a) A checker surface.
\qquad
b) The corresponding dessin d'enfant.
}
\label{fi:checker}
 \end{figure}
Consider a graph on this surface separating it into triangles. Let triangles be colored  black and white in the checker order, i.~e.,
neighbors of a black triangle are white and neighbors of  a white triangles are black. Let edges of  the graph be colored
 $3$ colors, denote them by $a$, $b$, $c$. Let these colors be arranged  clockwise  on a perimeter of each white triangles,
 let they be arranged anticlockwise on perimeters of black triangles.  We call such graphs on surfaces  by {\it checker surfaces}.
 Two surfaces are equivalent if they are isotopic.
So such a surface is a pure combinatorial object.

Denote by $\Kop$ the set of all checker surfaces, by $\Kop_n\subset \Kop$
we denote the
set of surfaces with $2n$ triangles ($n$ black amd $n$ white).

\sm

{\sc Labeled checker surfaces.}
Consider an element of $\Kop_n$. Assign pairwise different labels $1$, \dots, $n$ to 
black black triangles and pairwise different labels $1$, \dots, $n$ to white triangles. We call an object obtained in this way by a {\it labeled checker surface}.
Denote by $\Kop_n^\bullet$ the set of all such surfaces.

Fix an element $\Sigma\in\Kop_n^\bullet$. 
For each color $a$, $b$, $c$ we define an element $g_a$, $g_b$, $g_c\in S_n$ in the following
way. For each $k=1$, 2, \dots, $n$ we find a label $k$ on a white triangles. This  triangle has a unique edge of color $a$, this edge is
contained in a black triangle with some label, say $l$. Then $g_a$ sends $k$ to $l$. In the same way we define
permutations $g_b$ and $g_c$ (they are reflections through $b$-edges and $c$ edges respectively).
This determines a {\it one-to-one correspondence between the set {\rm $\Kop_n^\bullet$} and the group 
$S_n\times S_n \times S_n$.}

\begin{figure}
 $$\epsfbox{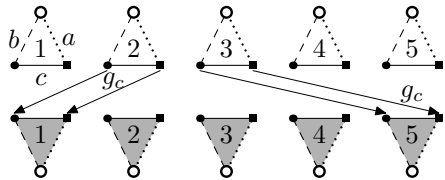}
 $$
\caption{Gluing of $\Sigma(g_a,g_b,g_c)$ from triangles.}
\label{fi:glueing}
 \end{figure}

{\sc The inverse construction.}
 Fix elements $g_a$, $g_b$, $g_c\in S_n$.
Consider a collection of $n$ white triangles with  
sides colored  $a$, $b$, $c$ clockwise, see Fig.~\ref{fi:glueing}. Assign to these triangles labels $1$, \dots, $n$. Similarly, take a collection of
$n$ labeled
black triangles, whose sides are colored $a$, $b$, $c$ counterclockwise.
For   each $k\le n$ we glue (according the orientations) 
the  side of type $a$ of $k$-th white triangle with side of the type $a$ of $g_a(k)$-th  a black triangle. 
Repeating this operation for $g_b$ and $g_c$,
we get an element of $\Kop_n^\bullet$. Denote a checker surface obtained in this way by
$\Sigma(g_a,g_b,g_c)$.

\sm

{\sc Types of vertices.}
Notice that we have 3 types of vertices. We say that a vertex has a type $A$ if it is adjacent to edges of colors $b$, $c$
(equivalently, it is opposite to an edge of the type $a$);
a type $B$ if it is adjacent to edges of colors $a$, $c$; a type $C$ if colors of adjacent edges are $a$ and $b$.

Colors of edges containing a vertex interchange, so
valences of vertices of checker surfaces are even. We say that the {\it order $\ord(v)$ of a vertex} $v$
 is the half of its valence.  For $\Sigma(g_a,g_b,g_c)\in \Kop^\bullet_n$,
 cycles of $g_a^{-1} g_b$ are in one-to-one correspondence with vertices of the type $C$, etc.
 Below it is convenient to use a double terminology:  {\it vertices of the type $C$ = vertices of the type $(g_a,g_b)$,} etc.

\sm 

%Let $\Sigma\in \vphantom{\Kop_n}_m \Kop_n^\bullet$. Consider a polygon and a vertex $v$ that is contained in  edges of colors $\daleth_i$ and $\daleth_{i+1}$.
%Then all edges meeting in $v$ a  colored $\daleth_i$ or $\daleth_{i+1}$, and colors interchange. 
%So the set of vertices splits into $m$ types, namely $(\daleth_1,\daleth_2)$, \dots, $(\daleth_{m-1},\daleth_m)$, $(\daleth_m,\daleth_1)$.

%Moreover, the set of vertices of a type $(\daleth_i,\daleth_{i+1}$ is in a one-to-one correspondence with the set of cycles of the permutation
%$g_{i+1}^{-1}g_i$. Labels on white triangles containing a vertex $v$ in counterclockwise order are elements of the cicle corresponding to $v$. Below we also call $(\daleth_{i+1},\daleth_i)$-vertices
%by {\it $(g_{i+1}, g_i)$-vertices}. Also we can call edges of a color $\daleth_i$ by $g_i$-edges.

\sm 

{\bf\punct Remark. Dessins d'enfant.} The set  $\Kop$ is in one-to-one correspondence with 'dessins d'enfant' in the sense of Grothendieck. Recall that a {\it dessin d'enfant}
is a bipartite graph on an oriented compact two-dimensional surface%
\footnote{Usually such surfaces are assumed to be connected, we
	admit disconnected   surfaces.}
such that the complement of the graph is a union of disks%
\footnote{Vertices are colored two colors, and edges can connect only vertices of
different colors.}. According the Belyi theorem,  dessin d'enfants are in one-to-one correspondence
with pairs $(R,\Phi)$, where $R$ is a non-singular complex curve defined over the algebraic closure $\ov \Q$ of rational numbers $\Q$, and $\Phi$ is a holomorphic function 
from $R$ to the Riemann sphere $\ov \C=\C\cup\infty$, whose critical values are $0$, $1$, $\infty$. For  detailed discussions, applications, and references, see \cite{LZ}, \cite{Sha}. 

Consider a checker  surface.
Removing edges of colors $a$, $b$ we get a dessin. Vertices of the type $C$
are in one to-one correspondence with complementary disks.

\sm

{\bf \punct Dirichlet measures on simplices.} Fix $\ell>0$.
Consider the $(p-1)$-dimensional simplex $\Delta_p(\ell)\subset \R^p$ defined by 
$$
x_1+\dots+x_p=\ell,\qquad  \text{where $x_j\ge 0$.}
$$ 
Let $k_1$, \dots, $k_p$ be positive reals. % let  $n:=\sum k_j$.
A {\it Dirichlet distribution} $\Theta_p[k_1,\dots,k_p;\,\ell]$ is a (probabilistic) measure on $\R^p$
supported  by $\Delta_p(\ell)$ and given by
$$
\frac{\Gamma(k_1+\dots+k_p)\,\ell^{-p+1} }{\prod_{i=1}^p  \Gamma(k_i)} \prod_{i=1}^p\, x_i^{k_i-1} dx_1\dots dx_{p-1},
$$
 see, e.g.,\cite{Wil}, 7.7.
This is a general definition, actually, we need only integer $k_j$

\sm

{\sc Remark.} In this formula $dx_p$ is absent, but the coordinate
$x_p$ is not distinguished, all measures 
$dx_1\dots dx_{j-1}\,dx_{j+1}\dots dx_p$ coincide on $\Delta_p(\ell)$.
\hfill $\boxtimes$

\sm

If $p=1$, we define $\Theta_1[k;\ell]$ as the delta-measure on $\R^1$ supported by $x_1=\ell$.

More generally, let $m_1$, \dots, $m_q\ge 0$. % denote $N=\sum m_j$.
Let $I$ be set of all $j$ such that $m_j>0$, say $I=\{m_{\alpha_1},\dots, m_{\alpha_p}\}$.
Then  $\R^q$ splits as a product of the space $\R^{\# I}$ with coordinates $x_\beta$, where $\beta\in I$,
and the space $\R^{q-\# I}$ with coordinates $x_\gamma$, where $\gamma\not\in I$.
We define $\Theta_q[m_1,\dots,m_q;\,\ell]$ as a product of the measure
$\Theta_{\# I}\bigl[m_{\alpha_1}, \dots, m_{\alpha_p};\,\ell]$ on $\R^{\# I}$ and the atomic unit measure  on $\R^{q-\# I}$ supported by
0.

\sm 

We need Dirichlet distributions for the following reason. Let $n$ points 
be distributed
uniformly on the circle of length $\ell$. They split the circle into $n$ pieces. Let $y_1$, \dots, $y_n$
be their lengths, such collections are elements of the simplex  $\Delta_n(\ell)$ distributed according $\Theta_n(1,\dots,1;\ell)$,
i.e., by
\begin{equation}
(n-1)!\,\ell^{-n+1}\,dy_1\dots dy_{n-1},
\label{eq:dxdx}
\end{equation}
this expression  is symmetric with respect to permutations of $y_1$, \dots, $y_n$.

\begin{lemma} 
	\label{l:dirichlet}
	Fix $k_1$, \dots $k_p\in \N$, $\sum k_j=n$.
	Consider the map $\Delta_n(\ell)\to \Delta_p(\ell)$ given by
	\begin{equation}
		x_1=\sum_{j=1}^{k_1} y_j,\quad
	x_2=\sum_{j=k_1+1}^{k_1+k_2} y_j,\quad
	\dots,\quad
	x_p=\sum_{j=k_1+\dots+k_{p-1}+1}^{n}.
	\label{eq:zamena}
	\end{equation}
Then the pushforward of the measure {\rm(\ref{eq:dxdx})} is the distribution
 $\Theta_p[k_1,\dots,k_p;\ell]$.
%\begin{equation}
%\frac{(n-1)!\ell^{-n+1}}{\prod (k_j-1)_j!} \prod x_j^{k_j-1} \,dx_1\dots dx_{p-1}
%\label{eq:push}
%\end{equation}
\end{lemma}

%{\sc Proof.}
%We pass to coordinates 
%\begin{multline}
%y_1, \dots, y_{k_1-1},  x_1,\,\, y_{k_1+1}, \dots, y_{k_1+k_2-1}, x_2,\,\,
%\\
% y_{k_1+k_2+1}, \dots,  y_{k_1+k_2+k_3-1}, x_3\,\,
%\dots,
%\label{eq:xyx}
%\end{multline}
%for this purpose we express variables $y_{k_1}$, $y_{k_1+k_2}$, \dots, $y_n=y_{k_1+\dots+k_p}$
%in terms of variables (\ref{eq:xyx}). The Jacobian of this transformation is 1.
%We get a body defined by the collection of inequalities:
%$$
%y_i\ge 0, \qquad\sum_{j=k_1+\dots+k_i+1}^{k_1+\dots+k_{i+1}} y_j\le x_i, \qquad \sum x_i=\ell.
%$$
% Now our map becomes a projection to a coordinate subspace
%$\{x_j\}$.
%By the Fubbini theorem, the pushforward of the measure is
%$$\ell^{-n+1} (n-1)!\, \Bigl\{\text{volume of a fiber}\Bigr\}\,dx_1\dots dx_{p-1}.
%$$
%A fiber over $(y_1,\dots,y_p)$ is the product of simplices
%$$ y_\alpha\ge 0\qquad(k_1+\dots+k_{i-1}<i<k_1+\dots+k_{i}),
%\qquad
%\sum_{\alpha k_1+\dots+k_{i-1}+1}^{k_1+\dots+k_{i}-1}
%y_\alpha<x_i.
%$$
%The volume of the product of these simplices is
%$\prod { x_j^{k_j-1}}/{(k_j-1)!}$.
%\hfill $\square$

This is a special case of the aggregation property of the Dirichlet distributions, see, e.g.,
\cite{Wil}, Subsect.~7.7.5.
\hfill $\square$

\sm 

{\bf \punct Laplace transforms of the Dirichlet distributions.}
Let $\Psi$ be a measure supported by $\R^p_+$. Its Laplace
transform  is 
$$
\cL \Psi(u_1,\dots,u_p):=\int_{\R^p_+} e^{-(u_1x_1+\dots+u_px_p)}\,d\Psi(x_1,\dots,x_p),
$$
where $\Re u_j\ge 0$.

\begin{proposition}
\label{pr:laplace}
	Let $k_j>0$.
	The Laplace transform of the Dirichlet distribution $\Theta_p[k_1,\dots,k_p;a]$
	is given by
	\begin{equation}
	\cL \Theta_p[k_1,\dots,k_p;a](u_1,\dots,u_p)=
	\frac{\Gamma\bigl(\sum k_j\bigr)\,a^{-p+1}}{2\pi i}
	 \int_{-i\infty}^{i\infty} e^{az} \prod_{j=1}^p (z+u_j)^{-k_j} \,dz.
	\label{eq:laplace1}
	\end{equation}
	If $k_j\in \N$, then this expession equals
	\begin{equation}
\frac{\bigl(\sum k_j-1\bigr)!\,a^{-p+1}}
{\prod_j (k_j-1)!}
  \prod_{j=1}^p\Bigl(-\frac\partial{\partial u_j}
\Bigr)^{k_j-1}
\Bigl[\sum_m \frac{e^{-au_m}}{\prod_{j:\,j\ne m}(u_j-u_m)}\Bigr].
\label{eq:laplace2}
	\end{equation}
\end{proposition}

Below we need  (\ref{eq:laplace2}). 
The statement (\ref{eq:laplace1}) was obtained in Phillips \cite{Phi},
we present proofs of both formulas.

\sm 

{\sc Proof.} Consider functions $f_1(x_1)$, \dots, $f_p(x_p)$
on $\R_+$ and an integral
$$
\phi_a:=\int_{\Delta_p(a)} f_1(x_1)\dots f_p(x_p)\,dx_1\dots dx_{p-1}.
$$
Consider its Laplace transform ($\Re z\ge 0$)
$$
\int_0^\infty e^{-za}\phi_a \,da=\int_{\R_+^p} e^{-z\sum x_j} \prod f_j(x_j) \,dx_1\dots dx_{p}
=\prod_{j=1}^p \int_0^\infty e^{-z x} f_j(x)\,dx.
$$
Applying the inversion formula for the Laplace transform, we get (cf. \cite{PBM}, formula (3.3.4.1))
\begin{equation}
\phi_a=\frac 1{2\pi i}\int_{-i\infty}^{i\infty} e^{za}
\Bigl(\prod_{j=1}^p \int_0^\infty e^{-z x} f_j(x)\,dx  \Bigr)\, dz.
\label{eq:phi-a}
\end{equation}

If $f_j(x_j)=x_j^{k_j-1} e^{-u_j x_j}$, then $\phi_a$ is (upto a constant factor)
the Laplace transform of $\Theta_p[k_1,\dots,k_p]$. 
We have
$$
 \int_0^\infty e^{-z x} f_j(x)\,dx =  \int_0^\infty
 x^{k_j-1} e^{-(z+u_j)x} \,dx=\frac{\Gamma(k_j)}{(z+u_j)^{k_j}}.
$$
Applying (\ref{eq:phi-a}), we come to the first statement  (\ref{eq:laplace1})
of the proposition. For general $k_j$
the integral in the right-hand side of (\ref{eq:laplace1}) 
 is a kind of multivariate confluent 
hypergeometric function. 
For integer $k_j$ this integral can be evaluated by residues.

Now set $k_j=1$ in (\ref{eq:laplace1}). Then
$$
\frac{1}{\prod_{j=1}^p (z+u_j)}=
\sum_m \frac1{\prod_{j:\,j\ne m}(u_j-u_m)} \cdot\frac{1}{z+u_m}.
$$
Keeping in mind the integral
$$
\frac{1}{2\pi i}\,  \mathrm{v.p.}\int_{-i\infty}^{i\infty} \frac{e^{az}\,dz}{z+u_m}=
e^{-au_m},
$$
we come to the identity
\begin{multline*}
\cL\Theta_p[1,\dots,1;a](u_1,\dots,u_p)=\\=
\int\limits_{\Delta_p(a)} \prod_{j=1}^p e^{-u_j x_j}\,dx_1\dots dx_{p-1}
= \sum_m \frac{e^{-au_m}}{\prod_{j:\,j\ne m}(u_j-u_m)}. 
\end{multline*}
A multiplication of a distribution
by $x_j$ implies a differentiation 
 of its Laplace transform by $u_j$. This leads
to (\ref{eq:laplace2}).
\hfill $\square$

\sm

{\bf \punct The action of $\cS[0,0]$ on $\frS_0^z$.%
\label{ss:0-0}}
Consider an element 
$$\{\bigcirc [k_j]\}_j\in \cS[0,0],$$
%$\frg[k_1,\dots,k_p]$ of $\cS(0,0)=\K\setminus \SS/\K$
 i.e.,  a union of cycles of lengths $k_1$, \dots, $k_p\ge 2$. To write a formula for the corresponding polymorphism, we need some notation.
 
 \sm 

Fix a collection of tables $\{U_\omega\}\in\frS^z_0$, let $\ell_\omega$ be their lengths, recall that they are pairwise
different almost sure.

%Denote it by $C(k_1,\dots, k_p)=C(\{k_i\})$.

For a given $\{U_\nu\}_\nu\in\frS^z_0$ we intend to write the spreaded image of $\{U_\nu\}_\nu$, i.e., the corresponding conditional measure on the space $\R^\circ \times \frS^z_0$.
This measure is supported by points of the type
$z^q\times \{V_\kappa \}_\kappa$, where a set $\{V_\kappa \}_\kappa$
coincides with $\{U_\nu\}_\nu$ up to a finite number of tables. It is more convenient to speak about collections
$\{\ell_{\nu}\}$ of lengths,
which are point configurations on the segment $[0,1]$. We introduce  the operation
\begin{equation*}
 \Replace\Bigl[\{\ell_{\omega_j} \}_j \longrightarrow \Lambda \Bigr] \{\ell_{\nu}\}
\end{equation*}
on such configurations, which replaces a finite subcollection $\{\ell_{\omega_j} \}_j $ of the collection $\{\ell_{\nu}\}$
by a random collection $\subset \R_+^p$ with a given distribution $\Lambda$.
More generally, we will use such notation for similar transformations of lists of other types.

\sm 

{\sc Notation.} a) For $\Sigma\in \Kop$ denote by $A_\alpha$ its $A$-vertices, by $B_\beta$ its $B$-vertices,
by $C_\gamma$ its $C$-vertices. Denote by $\#\{B_\beta\}$ (resp. $\#\{C_\gamma\}$) the number 
of $B$-vertices (resp. $C$-vertices).

\sm 

b) For $\Sigma\in \Kop$ and $\{ U_\nu  \}_{\nu\in \Omega}\in \frS_0^z$
we define a {\it framing} of $\Sigma\in \Kop$ as an injective map $\omega$ from the set of $B$-vertices
of $\Sigma$ to the set  of tables. In particular, for any $B$-vertex $B_\beta$ we assign a number $\ell_{\omega(B_\beta)}$. Denote by $\Fr(\Sigma,\{\ell_\nu\})$ the set of all framings of $\Sigma$.

\sm 

c) For $\Sigma\in \Kop$ we denote by $\Aut_B(\Sigma)$ the group of all automorphisms
of $\Sigma$ fixing all $B$-vertices%
\footnote{This group is poor, for a connected surface it is cyclic (and usually trivial),
	for a disconnected surface it is a product of cyclic groups.}.

\sm 

d) Let $\{k_j\}$ be as above. By $\Gamma[\{k_j\}]$ denote the set of all
$\Sigma\in \Kop_n$ whose $A$-vertices have orders $\{k_j\}$.

\sm 

e) Denote by
$\iota_m[\{k_j\}]$ the number of entries of $m$ to a list $\{k_j\}$.

%Denote by $\text{\Asterisk}_{j}$ the multiple convolution of measures on $\R^q$. By $x\Join y$ we denote a
%concatenation of two lists,
%$$
%(x_1,\dots,x_p)\Join (y_1,y_2,\dots):=(x_1,\dots,x_p,y_1,y_2,\dots).
%$$

\begin{theorem}
	\label{th:2}
	Fix an element $\{\bigcirc [k_j]\}_j\in \cS[0,0]$
	and a restaurant $\sqcup U_\nu$ with tables of lengths $\{\ell_\nu\}$.
%	, which is a product of cycles of lengths $\{k_j\}$.
% Let $\Gamma[k_1,\dots,k_p]$ be the set of all framed  checker surfaces $(\Sigma,\omega)$, whose $A$-vertices  have 
% orders $k_1$, \dots, $k_p$. 
% Let $B_\beta$ be  $B$-vertices of $\Sigma$ and  $C_\gamma$ be $C$-vertices. 
 For $\Sigma\in\Gamma[\{k_j\}]$ denote by  $m_{\beta\gamma}$ the number
 of edges connecting $B_\beta$ and $C_\gamma$.
 % Denote by $G(\Sigma, \omega)$ the group of automorphisms of $\Sigma$, preserving the framing,
 %i.e., fixing all $B$-vertices.%
 %\footnote{Nontriviality of this group is a relatively rare case. Such an automorphism preserves connected components; for a given connected
 %component it must induce a rotation about each vertex $B_j$.}
 %$B_i$.
 Then 
 the polymorphism
 $\Xi_{0,0}\bigl(\{\bigcirc [k_j]\}\bigr)$ sends the point $\{\ell_\nu \}$
 to the following measure 
 %$\xi_\ell{\{\ell(\omega)\}}$
  on $\R^\circ\times \frS_0^z$:  % given by the formula
 \begin{multline}
 \prod_j k_j \prod_{m\in \N} \iota_m[\{k_j\}]!
 \sum_{\Sigma\in \Gamma[\{k_j\}],\,\omega\in \Fr(\Sigma,\{\ell_\nu\})}
 \frac{1}  {\#\Aut_B(\Sigma)}
 \cdot \prod_{B_\beta} \frac{\ell_{\omega(B_\beta)}^{\ord(B_\beta)}} {(\ord(B_\beta)-1)!}
 \times\\\times
% \Add\Bigl[ \delta_{z^{r-q}}\times
\biggl(\delta_{\R^\circ} [[z^{ \#\{C_\gamma\} \,-
	\, \#\{B_\beta\}}]]\,\dot\times
\\ \, \dot\times\,
\Replace\Bigl[\bigl\{\ell_{\omega(B_\beta)}\bigr\}_\beta\longrightarrow
\Bigl(  \text{\footnotesize \AsteriskCenterOpen}_\beta\Theta_{\#\{C_\gamma\}}\bigl[\{m_{\beta\gamma}\}_\gamma;\,\ell_{\omega(B_\beta)}\bigr]\Bigr)
\Bigr]
\bigl\{\ell_{\nu}\bigr\}\biggr),
\label{eq:long1}
 %\coAsterisk
 \end{multline}
where the symbol \,\, {\footnotesize \AsteriskCenterOpen}  \,\, denotes a convolution of  Dirichlet  distributions on $\R^{\#\{C_\gamma\}}$.
\end{theorem}

Notice that a convolution in a summand of the formula is a  measure on the simplex
$\Delta_{\#(C_\gamma)}\bigl(\sum \ell_{\omega(B_\beta)}\bigr)$, densities of such measures
are piecewise polynomial.  In any case, by Proposition \ref{pr:laplace},
 Laplace transforms of such measures are  explicit elementary functions.
 
 \sm 

Theorem is proved in Subsect. \ref{ss:product-cycles}--\ref{ss:proof-th:2}.

\sm 

\begin{figure}
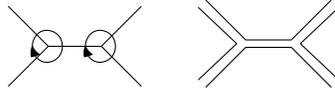

 	$$ \epsfbox{ribbon.4}  \qquad  \epsfbox{ribbon.5}
 	$$
 	
 	\caption{A ribbon graph and the corresponding band domain 
 	on a surface.}
 \label{fi:ribbon}
 \end{figure}

{\bf \punct A product of two permutations, which are determined in terms of disjoint cycles.%
\label{ss:product-cycles}}
Recall that a {\it ribbon graph} is a graph with fixed cyclic orders of edges at each vertex, see Fig.~\ref{fi:ribbon}. Such a graph can be regarded as an oriented two-dimensional surface consisting of 
bands along edges. Boundary of this surface is a union of circles,
 gluing a disk to each component of the boundary we get a compact
 oriented surface. 
 
 \sm 
 
\begin{figure}
	$$\epsfbox{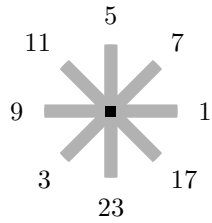}
	$$
	\caption{A chamomile.}
	\label{fi:chamomile}
\end{figure}

\begin{figure}
	$$\epsfbox{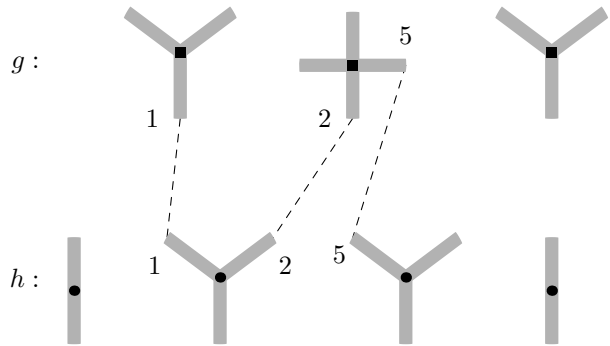}  
	$$
	\caption{Gluing of chamomiles.}
	\label{fi:cha-glue}
\end{figure}

{\sc A construction of a bipartite ribbon graph by two pemutations.}
 Let $g$, $h\in S_n$ be represented as products
of disjoint cycles
$$g=(\sigma_1)(\sigma_2)\dots , \qquad
h=(\tau_1)(\tau_2)\dots
$$
For each cycle $(\sigma_\alpha)$ we draw a '{\it chamomile}'
$A_\alpha$ as on Fig.~\ref{fi:chamomile}
and enumerate   petals of  $A_\alpha$ by elements of the cycle $(\sigma_\alpha)$
clockwise.  For each cycle $(\tau_\beta)$ of $h$ we draw a similar  chamomile $B_\beta$. Gluing petals of $g$-chamomiles and petals of $h$-chamomiles  according to the enumeration of petals and  orientations (see Fig.~\ref{fi:cha-glue}), we  get a bipartite ribbon graph,
say  $\Gr(g,h)$, whose list of vertices is $\{A_\alpha \}$, $\{B_\beta\}$ and edges
are enumerated by 1, 2, \dots, $n$.

\begin{figure}
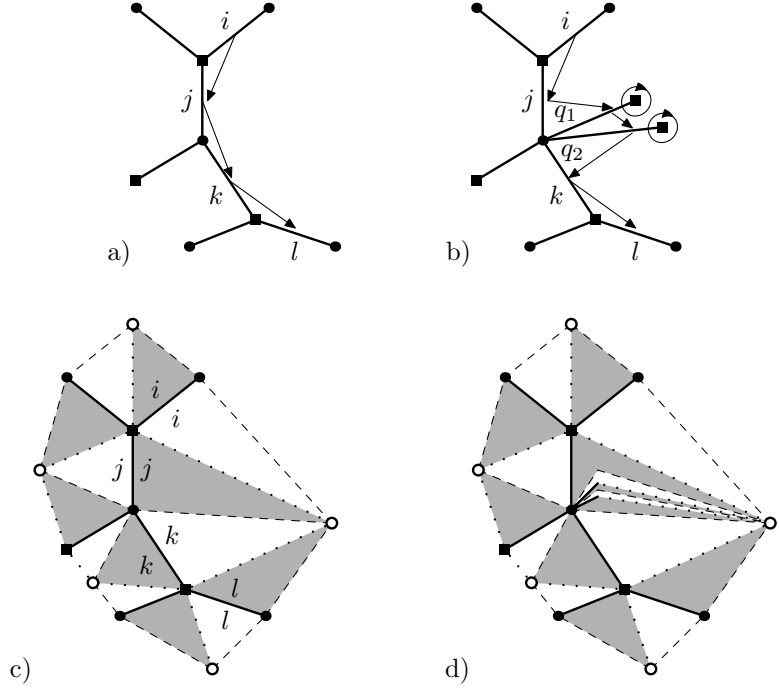

	$${\rm a)}\epsfbox{ribbon.8}\qquad\qquad
	{\rm b)} \epsfbox{ribbon.9}  
	$$

$${\rm c)}\epsfbox{ribbon.10}\qquad\qquad {\rm d)}  \epsfbox{ribbon.11} 
$$

\caption{a) A product of two permutations $g$, $h$ in terms of ribbon graphs. b) An insertion of two elements $q_1$, $q_2$ to a cycle of $h$ between $j$ and $k$. c) A labeled checker surface corresponding to elements $g$, $h^{-1}$, $1\in S_n$.
	 d) Insertion of two elements to a cycle of $h$ on the language of checker surfaces.}
\label{fi:insertion}
\end{figure}

Each component of the boundary of the ribbon graph $\Gr(g,h)$ can be regarded as a polygonal path
with an even number of sides. Passing it counterclockwise 
we observe edges of two types, with origins in $A$-vertices   and with origins in $B$-vertices. Labels on edges of the first type form cycles of the permutation $hg$, labels on edges of the second type form cycles of $gh$.
See Fig.~\ref{fi:insertion}.a.
 This correspondence arises at least to Goulden, Jackson \cite{GJ}.

\sm 

{\sc Insertions of elements to cycles and ribbon graphs.}
Such pictures are well-compatible with insertions of additional elements to cycles of $h$.
 Namely, take a larger group
$S_{n+N}$. Let $\wt g\in S_{n+N}$ be the trivial extension of $g$.
%\begin{equation}
%\wt g(m)=\begin{cases}
%g(m), \quad \text{for  $m\le n$}
%\\
%m,\quad \text{for  $m>n$}.
%\end{cases}
%\label{eq:wt-1}
%\end{equation}
Choose an element $h^{\circ N}\in S_{n+N}$ such that
\begin{equation}
\Upsilon^N_n(h^{\circ N})=h.
\label{eq:wt-2}
\end{equation}
Then the ribbon graph $\Gr(\wt g,h^{\circ N})$ can be easily obtained from the ribbon graph  $\Gr( g, h)$. 
Namely, an insertion of elements $q_1$, \dots, $q_r>n$ 
to a cycle $(\tau_{\beta})$  between $j$ and $k$,
$$
(\dots j\, k \dots)\mapsto (\dots j\,q_1\dots  q_r k \dots),
$$
means that we draw at the vertex $\beta$ edges labeled by $q_1$, \dots, $q_r$
between edges $j$ and $k$. Opposite vertices of these edges
have valence 1 (since $\wt g(q_\mu)=q_\mu$). See, Fig.~\ref{fi:insertion}.b.

\sm

{\bf\punct Checker surfaces and products of permutations.}
Reformulate this construction  in the terms of checker surfaces. For $g$, $h\in S_n$
we  take the surface $\Sigma(g,1,h^{-1})\in \Kop^\bullet_n$.
 Vertices of this surface  have types $(1,g)$, $(h^{-1},1)$,
$(g,h^{-1})$, they correspond to cycles of  permutations
$g$, $h$, $(hg)^{-1}$ respectively.

For each vertex of $(g,h^{-1})$-type,
labels in adjacent white triangles 
passing 
counterclockwise  give a  cycle of $hg$.

\sm 

Removing vertices of the type $(g,h^{-1})$ and adjacent edges
we get a ribbon graph described above, remaining edges
have type (1), labels on both sides of such an edge coincide and so these labels
 can be attributed to edges. We get the ribbon graph $\Gr(g,h^{-1})$.

\sm

{\sc Insertion of elements to cycles and checker surfaces.}
Operation of insertion of an element with number $>n$ to a cycle of
$h$ is shown on the  Fig.~\ref{fi:insertion}.

\begin{figure}
	$$\epsfbox{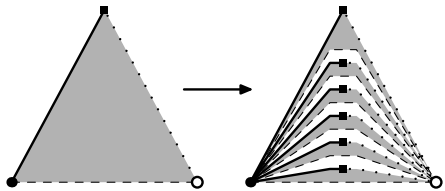}  
	$$
	\caption{To a pass from $\Sigma(g,1,h^{-1})$ to
		$\Sigma(\wt g,1, (h^{\circ N})^{-1})$. Additional drawing on former black triangles.}
	\label{fi:wt}
\end{figure}

So the surface $\Sigma(\wt g, 1, (h^{\circ N})^{-1})$
is obtained from the surface $\Sigma(g, 1,  h^{-1})$
by the following operation. We preserve all white triangles with labels, and black triangles are transformed
as it is shown on Fig.~\ref{fi:wt}.

\begin{figure}
	$$\epsfbox{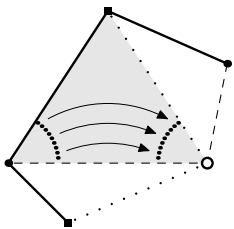} 
	$$
	\caption{Transformations of a piece of a table with guests.}
\label{fi:guests}
\end{figure}

{\bf \punct The right action of $S_n$ on the space of virtual permutations $\frS^z$.%
\label{ss:dense}} Let $g\in S_n$, $\fru\in \frS^z$. Denote $u=\Upsilon^\infty_n \fru$.
Consider the surface $\Sigma\bigl(g,1, u^{-1}\bigr)\in \Kop^\bullet_n$. 
Take a $B$-vertex, let $T$ be an adjacent black triangle.
Let the clockwise white neighbor of $T$ has label $j$
and so the counterclockwise white neighbor has label $u(j)$.
We assign to $T$ the interval 
$(\ov j, \ov{u(j)})$ of the corresponding occupied  table, denote this interval by $\cO(T)$, see Fig.~\ref{fi:guests}.

For each $C$-vertex consider a counterclockwise cyclic chain
$$j_1,\cO(T_{i_1}), j_2, \cO(T_{i_2}) , \dots$$
consisting of labels on adjacent white triangles and 
ordered sets $\cO(\cdot)$ on adjacent black triangles. Uniting them,
we get a table of $\fru g$.

\sm

 %This allows to give a geometric description a product of an element $g\in S_n$
%and a virtual permutation $\fru$. Namely, we take the permutation
%$u:=\Upsilon^\infty_n (\fru)$ and draw the labeled triangular
%surface corresponding to $g$ and $u^{-1}$. For each black triangle
%with label $k$ we take the table of $\fru$ containing $k$ and the interval $\ov k< \ov m< u(k)$ of this table. Next, we take a $(g,u^{-1})$-vertex
%of the surface and read anticlockwise labels on white triangles and semisegment
%on black triangles.

%This description is not simpler than the description of 'products' in Subsect. \ref{ss:virtual}. But this allows to write formulas for
%polymorphisms.

\sm

{\bf \punct The action of $S_n$ on $\frS_n^z$.} 
	For an element $g\in S_n$
	consider the polymorphism
	$$\Xi_{n,n}(g)=\frx[n]^\star g \frx[n]$$
	 of $\frS_n^z$.
	Consider a point $\fru$ of $\frS_n^z$.
It determines a permutation $u\in S_n$
and a collection of lengths $\{l_j\}$ between points $\ov j$ and $\ov{u(j)}$. 
	
	Consider the labeled checker surface $\Sigma=\Sigma(g, 1, u^{-1})\in \Kop_n^\bullet$. 
	For a black triangle of $\Sigma$ with label $j$ we assign the length of the  arc $[\ov j, \ov{u(j)}]$. 
	
	 \begin{lemma}
	 	\label{l:frT-n}
	 The spreaded image of a point $\fru\in \frS_n^z$ under the polymorphism $g\in S_n$
	 is a $\delta$-measure on $\R^\circ\times \frS^z_n$ supported by 
	 a point $z^s\times \frv$, where 
	 $$s=\#\{\text{\rm $(g^{-1},u) $-vertices}\}-\#\{\text{\rm $(u^{-1},1) $-vertices}\}$$
	 and  $\frv\in \frS_n^z$ is defined in the following way.
	  For each $(g,u^{-1})$-vertex
	 $C_\gamma$ of $\Sigma$ we draw a table taking labels $\ov \mu_1$, $\ov \mu_2$, \dots from white triangles passing counterclockwise and  length of arc $[\ov\mu_{i-1},\ov \mu_i]$ from black triangles. This gives us a collection
	 of tables containing all labels $\ov 1$, \dots, $\ov n$. Adding
	 tables of $\fru$ that do not contain labels we get $\frv$.
	 \end{lemma} 
 
 In particular, in this case we get a deterministic map $\frS_n^z\to\frS_n^z$,
 so the group $S_n$ acts  on the space $\frS^z_n$.
 
 \sm 
	
%{\sc Proof.} Consider an element of the projective limit
%$\frS$ lying over $u$, i.e., a sequence $u^{\circ N}\in S_{N+n}$, such that $\Upsilon^N_{N-1}(u^{\circ N})=u^{\circ(N-1)}$, $u^{\circ 0}=g$.
%Let us require 
%$$
%\Upsilon^\infty_N\Bigr(\lim_{N\to\infty} u^{(N)}\Bigl)=\fru
%$$	
%By (\ref{eq:limit}) this means the following. For 
%$i\le n$ consider the cycle $\sigma_{[i]}$
%containing $i$ and  denote by $p_N^i(h^{\circ N})$ the number of new elements $j\le n+N$, which were inserted to  $\sigma_{[i]}$
%between $i$ and $u(i)$. 
%Then for all $i$
%\begin{equation}
%\lim_{N\to\infty} p^N_i(h^{\circ N})=\bigl\{ \text{length of $[\ov i, \ov{u(i)}]$} \bigr \}
%\label{eq:limlim}
%\end{equation}
%We draw surfaces $\Sigma(\wt g,1, (u^{\circ N})^{-1})$
%and get cyclic structure of elements $u^{\circ N}\wt g$.
% Fig.~\ref{fi:wt} shows a result of insertion of elements to cycles, we have the same (non-ordered) collection 
 %$\{p^N_i(h^{\circ N})\}_i$ of numbers of inserted elements,
 %and therefore limits of ratios
 %are the same.
 %\hfill $\square$

The statement follows from the 
the construction of the previous subsection 
and formula (\ref{eq:limit}).

%\sm

%{\bf \punct The right action of $S_n$ on $\frS^z$.} Let $g\in S_n$, $\fru\in \frS$. Denote $u=\Upsilon^\infty_n \fru$ Consider the surface $\Sigma\bigl(g,1, u^{-1}\bigr)\in \Gamma^\bullet_n$. 
%Take a $B$-vertex, let $T$ be a black triangle.
%Let the clockwise white neighbor of $T$ has label $j$.
%We assign to $T$ the interval 
%$(\ov j, \ov{u(j)})$ of the corresponding table considered as 
%a countable ordered set, denote this set by $\cO(T)$.

%For each $C$-vertex consider a counterclockwise chain
%$$j_1,\cO(T_{i_1}), j_2, \cO(T_{i_2}) , \dots$$
%consisting of labels on adjacent white triangles and 
%ordered sets $I(\cdot)$ on black triangles. Uniting them,
%we get a table of $\fru g$.

\sm

{\bf \punct Proof of Theorem \ref{th:2}.%
\label{ss:proof-th:2}}
Denote $n:=\sum k_j$.
Choose an element $g\in S_n$ whose cycles have lengths 
$k_1$, \dots, $k_p$. Denote by $(\sigma_\alpha)$  disjoint cycles
of $g$. This choice is noncanonical and further considerations 
depend on $g$ until the last paragraph of the proof.
  We extend $g$ to an element $\wt g$ of $S_\infty$ in a trivial way. Then the double coset
$\K\cdot (1,\wt g)\cdot \K$ is $\{\bigcirc[k_j] \}$.
We decompose
$$\{\bigcirc[k_j] \}=(\lambda^n_0)^* g \lambda^n_0.$$

\sm 

Recall that the polymorphism $\Xi(\lambda^n_0):\frS^z_0\to \frS_z^n$
defined in Subsect. \ref{ss:embeddings} is a random uniform arrangement of
guests $\ov 1$, \dots, $\ov n$ in an empty restaurant.
This operation 
can be described in the following (more complicated) way.

 Fix a restaurant, let 
$\{\ell_\mu\}$ be lengths of tables (recall that they are pairwise distinct a.s).
%Denote by $\cB_\mu$ 

--- First, fix a collection $\{s_\mu \}$ of nonnegative integers such that 
$\sum s_\mu=n$. Assume that a probability of such a choice
is
$$
\frac{n!\, \prod \ell_{\mu}^{s_\mu}}{\prod s_\mu!}
$$
(products actually are finite, the sum equals 1).
Denote $\cB=\cB(\{s_\mu\})$ the set of all $\mu$, for
which $s_\mu>0$.

\sm

--- Second, for each $s_\mu>0$ we take a cyclically ordered set $V_\mu$ with $s_\mu$
elements. Choose a bijection $\{1,2, \dots, n \}\to \sqcup_{\mu\in \cB} V_\mu$
and assume that all variants have the same probability $1/n!$. So we get a collection $\{V_\mu \}_{\mu\in\cB}$
of  cyclically ordered sets consisting of positive integers $\le n$.
We consider each set up to cyclic ordering, so there are $\prod_{\mu\in \cB} s_\mu$ ways to obtain
the same collection. Thus, we get a permutation $u\in S_n$ defined as a product
of  cycles.% and  also a collection of lengths $\ell_{\mu}$. 

\sm 

--- Third, for each cycle $V_\mu$ of $u$ we insert arcs between elements of the cycle (guests), 
lengths of arcs are chosen according
the uniform probabilistic distribution 
on the simplex $\Delta_{s_\mu}(\ell_\mu)$  (i.e., the Dirichlet distributon $\Theta_{s_{\mu}}[1,\dots,1;\ell_\mu]$). 

\sm

So  we get an element $\fru$ of $\frS_n^z$ depending on $\{s_\mu\}$,  $u\in S_n$, 
and a point of product of simplices $\prod_{\mu\in\cB} \Delta_{s_\mu}(\ell_{\mu})$.

\sm 

Keeping in mind Lemma \ref{l:frT-n}, we can  apply $g\in S_n$ to a random element $\fru\in \frS^z_n$.

For each $u$ consider the surface $\Sigma(g,1,u^{-1})$. By definition, labels on both sides of any
 edge of the color $(1)$ coincide.
Vertices of types $A$, $B$, $C$ in formulation of Theorem \ref{th:2} correspond to types
$(1,g)$, $(u^{-1},1)$, $(g,u^{-1})$. Clockwise cycles of labels
on white triangles about a vertex $A_\alpha$ are cycles $(\sigma_\alpha)$
of $g$.
Each  surface is equipped with a framing $\omega$, namely $B$-vertices correspond to tables $\{U_\mu\}_{\mu\in \cB}$
 of the initial restaurant $\{U_\mu\}$
(comparatively to the set $\Gamma[\{k_j\}]$ in the formulation of Theorem \ref{th:2}, we add labels on triangles).
 Additionally,  we
assign a positive number $l(j)$ to each black triangle $T(j)$ with label $j$. Such numbers around a given $B$-vertex
 are distributed according the measure 
 $\Theta_{\ord(B_\beta)}\bigl[1,\dots, 1;\ell_{\omega(B_\beta)}\bigr]$.

New tables of a random restaurant $\Xi_{0,n}((\lambda^n_0)^*)  \Xi_{n,n}(g) \Xi_{n,0}(\lambda^n_0)  $  correspond to vertices $C_\gamma$, a length $L(\gamma)$ of such a table
is sum of real numbers $l(j)$  attributed to black triangles containing $C_\gamma$,
$$
L(\gamma)=\sum_{j:\, C_\gamma\in T(j)} l(j).
$$
We write this sum in the form
$$
L(\gamma)=\sum_{B_\beta} \,\,\sum_{j:\,\, B_\beta,C_\gamma\in T(j)} l(j).
$$
By
Lemma \ref{l:dirichlet}, for a fixed $\beta$ a distribution of vectors 
$$
\Bigl\{\sum_{j:\,\, B_\beta,C_\gamma\in T(j)} l(j)\Bigr\}_\gamma
$$
is the Dirichlet distribution 
$\Theta_{\#(C_\gamma)}[\{m_{\beta\gamma}\}_\gamma; \ell_{\omega(C_\beta)}]$.
Therefore the sums $\{L(\gamma)\}_\gamma$ are distributed as a convolution
$\text{\footnotesize\AsteriskCenterOpen}_\beta$ of such Dirichlet distributions.

Hence the desired conditional measure
is
\begin{multline*}
\sum_{\text{$\Sigma(g,1,u^{-1})$ with given $g$}} \quad
\sum_{\omega\in \Fr(\Sigma(g,1,u^{-1}),\{\ell_\nu\})}
\,\,
\prod_{B_\beta} \frac{\ell_{\omega(B_\beta)}^{\ord(B_\beta)}} {\ord(B_{\beta})-1)!}
 \times\\\times
 \Biggl( \delta_{\R^\circ} [[z^{ \#\{C_\gamma\}  \,-\, \#\{(B_\beta)\}} ]]\,
 \dot\times\\\dot\times\,
 \Replace\Bigl[\bigl\{\ell_{\omega(B_\beta)}\bigr\}\longrightarrow
\Bigl(  \text{\footnotesize \AsteriskCenterOpen}_\beta\Theta_{\#\{C_\gamma \}}\bigl[ \{m_{\beta\gamma}\}_\gamma; \ell_{\omega(B_\beta)}\bigr]\Bigr)\Bigr]
\bigl\{\ell_\nu\bigr\}
\biggr).
\end{multline*}

In fact, the first summation is taken over the subset in $\Kop_n^\bullet$
consisting of labeled surfaces such that: 

\sm 

1) for each $c$-edge labels on two-sides of the edge coincide;

\sm 

2) clockwise cycles of labels of white triangles   
about $A$-vertices coincide with the disjoint cycles $(\sigma_\alpha)$ of the
permutation $g$.

\sm 

Forgetting all labels we get an element of the set $\Gamma[\{k_j\}]$ defined 
in  Subsect. \ref{ss:0-0}. Fix a surface  $\Lambda\in\Gamma[\{k_j\}]$
with fixed framing (so we can think that $B$-vertices of $\Lambda$ are enumerated)
and evaluate the number of different surfaces of the type
$\Sigma(g,1,u^{-1})$ over $\Sigma$. 

\sm 

1) We must choose a one-to-one correspondence
between the set of $A$-vertices and the set of cycles $(\sigma_\alpha)$
of $g$ (an order of a vertex must coincide with an order of the corresponding cycle).  
This can be done by
$\prod_{m\ge 2} \iota_m\bigr[\{k_j\}\bigr]!$ ways.

\sm

2) After this, the collection of labels on white triangles containing a vertex
$A_\alpha$ consists of elements of the corresponding cycle
$(\sigma_\alpha)$ in the same cyclic order. The only freedom is rotation of the cycle,
this gives us the factor $\prod \ord(A_\alpha)=\prod k_\alpha$

\sm 

3) Labels on black triangles are uniquely determined by labels on white triangles.

\sm 

4) Two surfaces $\wt \Lambda_1$, $\wt \Lambda_2$ obtained by such arrangements  (including framings) can 
be equivalent. In this case we have an automorphism of the surface preserving the framing
and identifying two labelings. So we must divide the result by the order of the group
of such automorphisms.

\sm

We came to formula (\ref{eq:long1}).
\hfill $\square$

%Next, let us pass from $\Gamma^\bullet[\dots]$ to $\Gamma[\dots]$. Consider a framed surface
%$(\Sigma,\omega)$ of a given type without labels. For an arrangement of labels
%we must set one-to-one correspondence between set of cycles of $g$ and $A$-vertices
%of $\Sigma$ according orders. This can be done by
%$\prod_{m\in\N}\iota(\{k_j\})!$ times. Also, we have a freedom of rotation of labels about 
% each $A$-vertex. This gives a factor $\prod k_j$. On the other hand, some arrangements
% can be identified by an automorphism of $\Sigma$ preserving the framing. For this reason,
% we must divide a summand by the order of group of such automorphisms.
% This gives us formula (\ref{eq:long1}).
% \hfill $\square$
 
\section{Further statements}

\COUNTERS

{\bf\punct A more general case.\label{ss:left-trivial}%
} Next, let us consider a space $\frS_n^z$
and describe its polymorphism determined by a chip, whose
 left-hand half is trivial, see Fig.~\ref{fi:left-trivial}.
 \begin{figure}
  $$\epsfbox{chips-add.3}
  $$
  \caption{To Subsect. \ref{ss:left-trivial}. A chip $\frb$}
  \label{fi:left-trivial}
 \end{figure}

 \sm

{\sc Data defining a chip.}
Fix a permutation $\sigma\in S_n$, a collection
$\phi_1$, \dots, $\phi_n\in \Z_+$ and a collection of integers
$k_1$, \dots, $k_p\ge 2$. Define a chip 
$$\frb=\frb\bigl[\sigma, \{\phi_i\}, \{k_j\}\bigr]\in \cS(n,n)$$
 consisting of the following arcs:

\sm 

1) $j_l^+\downarrow j_l^-[0]$, where $j=1$, \dots, $n$;

\sm 

2) $j_r^+\downarrow \sigma(j)_r^-[\phi_j]$, where where $j=1$, \dots, $n$;

\sm 

3) $\bigcirc[k_i]$, where $i=1$, \dots, $p$.

\sm 

{\sc Data determining a half-empty table.}
Next, fix an element $\tau\in S_n$ represented as a product of disjoint cycles, a collection of positive reals $l_1$, \dots, $l_n$,
and a countable collection of positive reals $\ell_\nu$ such that $\sum_i l_i+\sum \ell_\nu=1$.
These data define  a point 
$$\frt=\frt[\tau,\{l_i\}, \{\ell_\nu \}]\in \frT_n^z.$$
 Namely, for $i$-th cycle of $\tau$
we create a table $U_i$ containing points of the cycle as guests  (in the same cyclic order), and claim
that distances between $j$ and $\tau(j)$ are $l_j$. Also, we add label-free tables
of lengths $\ell_\omega$.

\sm

Now we intend to write the spreaded image of the half-empty restaurant $\frt$ under 
the polymorphism $\frb$.
The result is similar to Theorem \ref{th:2}, but it must include  additional
combinatorial  data.

\sm 

{\sc The set of summation.}
Define the following set
$\Gamma[\frb;\frt]$
consisting of  checker surfaces equipped with following additional 
structures: 

\sm

1) An  injective map $F$ ({\it labeling})
$$
F:\,\bigl\{1,2,\dots, n \bigr\}
\rightarrow\,
\bigr\{\text{set of white triangles}\bigr\}.
$$

\sm

\begin{figure}
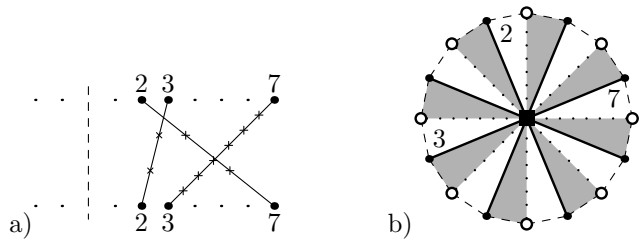

	$${\rm a)}\epsfbox{chips-add.4}\qquad\qquad
	{\rm b)} \epsfbox{ribbon.14}
	$$
	\caption{a) A cycle of $\sigma$ (on the figure we have a cycle (2\,7\,3)) and  corresponding arcs of
		a chip. \, b) The $A$-vertex corresponding to this cycle $\sigma$.}
	\label{fi:A-vertex}
\end{figure}

2) A bijection between 
\begin{multline*}
\bigl\{\text{set of $A$-vertices of $\Sigma$}\bigr\}
\longleftrightarrow
\\
\longleftrightarrow
\bigl\{
\text{set
of cycles of $\sigma$}\bigr\}\cup \bigl\{\text{set of  cycles $\bigcirc[k_i]$ of the chip $\frb$}\bigr\}.
\end{multline*}
This map must be compatible with the labeling $F$. Namely, if a vertex $A_\alpha$ corresponds to a cycle $\bigcirc[k_i]$ of the chip, then $\ord A_\alpha=k_i$
and white triangles containing $A_\alpha$ are label-less.
If a vertex $A_\alpha$ corresponds to a cycle $(\nu_1\,\nu_2\, \dots\, \nu_\rho)$ 
of $\sigma$, then labeled white triangles containing $A$ are precisely   $F(\nu_1)$, \dots, $F(\nu_\rho)$ 
in this cyclic order, and for each $\nu_i$ there are precisely $\phi_i$ label-less
white triangles between  $F(\nu_i)$ and $F(\nu_{i+1})$, see Fig.~\ref{fi:A-vertex}.

\begin{figure}
	$$\epsfbox{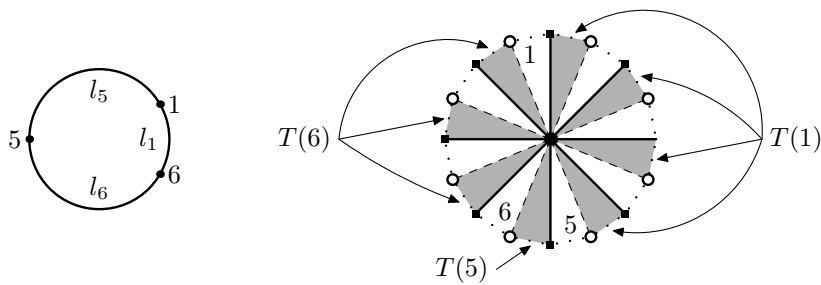}
	$$
	\caption{A table, the corresponding $B$-vertex, and sets $T(i)$.}
	\label{fi:Ti}
\end{figure}

\sm 

3) An injective map $\omega$ (framing) 
$$
\omega:\,\bigl\{\text{set of $B$-vertices 
of $\Sigma$}\bigr\} \rightarrow\,  \bigl\{\text{set of tables of $\frt$}\}.
$$
Moreover, labels on white
triangles containing $B_\beta$ must coincide with guests on the corresponding table
in reversed cyclic order, see Fig.~\ref{fi:Ti}.

\sm 

%This $\Gamma[\frb;\frt]$  replaces the summation index $\Gamma(\{k_i\})$ in Theorem \ref{th:2}.
%Next, we explain how to write summands.

{\sc Some notation.}
1)  Vertices $B_\beta$ split into two types. Denote by $\Phi$ the set of all $B_\beta$
such that all white triangles containing $B_\beta$ are label-less.
For $B_\beta\in \Phi$ denote by $T(B_\beta)$ the set of {\it black} triangles containing $B_\beta$.

For $B_\beta\notin \Phi$ consider the counterclockwise cycle of adjacent white triangles,
let $\kappa_1$, \dots, $\kappa_s$ be labels on labeled white triangles. For each adjacent pair
$(\kappa_i,\kappa_{i+1})$ denote by $T(\kappa_i)$ the set of {\it black} triangles
between them. Notice that a $B$-vertex  is determined by the number $\kappa_i$ and we can omit $\beta$ from the notation, see Fig.~\ref{fi:Ti}.

\sm 

2) Similarly, denote by $\Psi$ the set of $C$-vertices, for which all adjacent triangles
are label-less. For $C_\gamma\in \Psi$ denote by $T(C_\gamma)$ the set of {\it black}
triangles containing $C$. For $C\notin \Psi$ we write counterclockwise labels on adjacent
white triangles $(\mu_1,\mu_2,\dots)$. Denote by $R(\mu_j)$ the set of {\it black} triangles
between $\mu_j$ and $\mu_{j+1}$.

\sm 

3) For each vertex of $C$ we consider the counterclockwise cycle of adjacent labeled triangles. Uniting
such cycles we get an element $\rho$ of $S_n$.

\sm

4) Denote by $\Aut_B(\Sigma,F)$ the group of all automorphisms of $\Sigma$ fixing all $B$-vertices
and all labeled triangles. Such an automorphism   automatically is trivial on each
connected component of $\Sigma$ containing a labeled triangle. 

\sm 

{\sc Remark.}
In Theorem \ref{th:2} parameters of Dirichlet measures $\Theta[\dots]$ are expressed in the terms of the matrix $m_{\beta\gamma}$ whose 
elements are number of black triangles with vertices $B_\beta$ and $C_\gamma$.
Rows of the matrix are enumerated by $B$-vertices and columns by $C$-vertices. In our case, 
rows of a similar matrix are enumerated by chains $T(i)$ between two labeled triangles at $B$-vertices $\notin \Phi$ and by vertices $B_\beta\in\Phi$.
 Columns are enumerated by chains $R(j)$ between labeled triangles at $C$-vertices $\notin \Psi$ and by vertices  $C_\gamma\in \Psi$.  So analog of the matrix $m_{\beta\gamma}$ contains elements of types
$$
\begin{matrix}
\#(T(i)\cap R(j))& \#(T(i)\cap R(C_\gamma))\\
 \#(T(B_\beta)\cap R(j)) & \#(T(B_\beta)\cap R(C_\gamma))
\end{matrix}
$$
For this reason an expression  for convolutions now is longer, but in fact the structure
of the formula is the same.
\hfill $\boxtimes$

\begin{theorem}
\label{th:3}
	The spreaded image of a point ${\frt}$ under a polymorphism $\frb$ is the measure on $\R^\circ\times \frS_n$
	 given by the formula
\begin{multline}
\prod k_i \prod \iota_m[\{k_j\}] \sum_{(\Sigma,\omega,F)\in\Gamma[\frb,\frt]}
\frac 1{\#G(\Sigma,F)}
%(n+\sum k_i+\sum\delta_j)! \,
 \prod_{\nu=1}^n \ell_\nu^{\# T(\nu)}
	\prod_{B_\beta\in \Phi}\frac{\ell_{\omega(B_\beta)}^{\ord(B_\beta)}} {(\ord (B_\beta)-1)!}
\times\\\times
\biggl(\delta_{\R^\circ}\bigl[\bigr[z^{\#\{C_\gamma\}- \#\{B_\beta\}}\bigr]\bigr]
\dot\times 
 \\ \dot\times
\Replace\Bigl[(\tau;\{l_i\};\{\ell_{\omega(B_\beta)}\}_{\beta\in\Phi}  )
\longrightarrow
\\
\Bigl(\rho;
 \underset{i\in I_n}{\text{\footnotesize\AsteriskCenterOpen}}
 \Theta_{n+\# \Psi}\bigl[ \{\#(T(i)\cap R(j))\}_{j\in I_n},
 \{\#(T(i)\cap R(C_\gamma) )\}_{C_\gamma\in \Psi};l_i
 \bigr]
 \Bigr)
 \text{\footnotesize\AsteriskCenterOpen}
\\
\Bigl(\underset{B_\beta\in\Phi}{\text{\footnotesize\AsteriskCenterOpen}}
%\SVER_{B_\beta\in\Phi}
\Theta_{n+\# \Psi}\bigl[\{\#(T(B_\beta)\cap R(j)) \}_{j\in I_n},
 \{\#(T(B_\beta)\cap R(C_\gamma)) \}_{C_\gamma\in \Psi}  ;\ell_{\omega(B_\beta)}\bigr]\Bigr)\Bigr]
 \\ 
 \bigl(\tau;\{l_i\};\{\ell_\nu\} \bigr)
 \biggr).
\end{multline}
\end{theorem}

For a proof we
set $N=n+\sum \phi_i+ \sum k_j$ and
take an element $g\in S_N$ such that
$$
\frb=(\lambda^N_n)^* g\lambda_n^N.
$$
Further arguments repeat the
proof of Theorem \ref{th:2} and we omit them.

\sm 

{\bf \punct Action of the center of $\cS(\infty,\infty)$
	on $\frS^z$.} Let $k_j\ge 2$.
Now consider the  element   $\mathbf{1}[\{k_j\}]\in\cS(\infty,\infty)$ consisting of

\sm 

1) arcs $i_l^+\downarrow i_l^-[0]$,
$i_r^+\downarrow i_r^-[0]$ for $i\in \N$
(this collection of arcs defines the unit in $\cS(\infty,\infty)$;

\sm 

2) cycles $\bigcirc[k_j]$. 

\sm

Such elements form the center of the semigroup $\cS(\infty,\infty)$.

Let $\fru\in \frS^z$, it is a collection $\{\ov {U_\nu}\}$
of tables equipped with guests.
Consider the set $\Gamma^\diamond[\{k_i\}]$ of framed  checker surfaces whose $A$-verices have orders $k_j$.
Fix $(\Sigma,\omega)\in \Gamma^\diamond[\{k_i\}]$.
For each vertex $B_\beta$ denote by $Y(B_\beta)$ the set of all
white triangles containing $B_\beta$. For each $B_\beta$
consider a random map $\eta_\beta:Y(B_\beta)\to U_{\omega(B_\beta)}$ reversing the cyclic order. For each $\beta$ we equip the set of all $\eta_\beta$ with the uniform
measure $d\eta_{\beta}$, assuming that the total measure is $\ell_{\omega(B)}^{\ord B_\beta}/(\# B_\beta-1)!$. Let $Q$ be a black triangle containing
$B_\beta$. 
Let $T$ be its white predecessor according to the cyclic order, $T'$ be the white follower. We assign 
to $Q$ the arc $\cJ(Q):=[\eta_\beta(T), \eta_\beta(T')]$ of the table $U_{\omega(B_\beta)}$ together with guests. 
Notice that the ends of such arcs are not guests a.s.
Define a new collection  $\{\zeta_\gamma\}=\{\zeta_\gamma(\{\eta_\beta\})\}$
of occupied tables
in the following way: take a vertex $C_\gamma$, write 
all black triangles $Q_{s_1}$, $Q_{s_2}$, \dots containing $C_\gamma$ in the counterclockwise order.
Then $\zeta_\gamma$ is a cyclically ordered set
$$
\zeta_\beta:=\cJ(Q_{s_1})\sqcup\cJ(Q_{s_2})\sqcup\dots
$$
equipped with a natural length. We also have a pushforward
$\zeta_*(d\eta)$
of the measure $d\eta=d\eta_1 d\eta_2\dots$, this is a new measure on a collection of occupied tables.

\begin{theorem}
	The image of a point $\fru\in \frS^z$ under the polymorphism $\mathbf{1}[\{k_j\}]$
	is
	\begin{multline*}
	\prod_j k_j! \prod_{m\ge 2} \iota_m[\{k_j\}]	\sum_{(\Sigma,\omega)\in \Gamma^\diamond[\{k_j\}]}
	 %\ell_\nu^{\# T(\nu)-1}
	%\prod_{B_\beta\in \Phi}\ell(\omega(B_\beta))^{\ord(B_\beta)} \ord (B_\beta)
	\\
     \Bigl(	\delta_{\R^\circ}\bigl[\!\bigl[ z^{\#\{C_\gamma\}- \#\{B_\gamma\}}\bigr]\!\bigr]
	\dot\times 
	\Replace\Bigl[
	\{\ov{U_{\omega(B_\beta)}}\}\longrightarrow 
	\{\zeta_*(d\eta)\}
	\Bigr]\{\ov{U_\nu}\}\Bigr)
	\end{multline*} 
\end{theorem}

To evaluate this image we apply Remark after Proposition \ref{pr:closure}. We take the second copy  $\un \N$ of $\N$ consisting of points $\underline 1$, $\underline 2$, \dots
Let regard $\ov S_\infty$ as the group $\ov S(\N\sqcup \un \N)$ of  all permutations
of $\N\sqcup \un \N$ and $\SS$ as a subgroup in 
 $\ov S(\N\sqcup \un\N)\times \ov S(\N\sqcup \un\N)$.
 Now we can realize $\cS(\infty,\infty)$ as the semigroup
 of double cosets $\ov S(\underline \N)\backslash \SS/\ov S(\underline \N)$ 
 
 Let $\sum k_j=n$.
 Take $g\in \SS$ that is trivial on on the left copy of
 $\N\sqcup\underline\N$, on right copy of $\N$ and on the set
 $\{\un{n+1}, \un{n+2}, \dots\}$. Let also $g$ induces a permutation of the set $\{\un 1,\dots, \un n\}$ with cycles of lengths $k_j$.
 
 Further considerations as the same as above. 
 We take an analog of $\lambda^n_0$, it corresponds
 to the chip $\un 1^-_l[1/2]\smallfrown \un 1^+_r[1/2]$, 
 \dots, $\un n^-_l[1/2]\smallfrown \un n^+_r[1/2]$
 and other arcs are strictly vertical lines of weight 0. It corresponds to putting additional
 guests $\un{\ov 1}$,\dots, $\un{\ov n}$
 to the occupied restaurant $\fru$.
 
 We fall to
 the situation of Subsect. \ref{ss:dense}, where the set
 $\N$ is replaced by $\N\sqcup\{\un{\ov 1},\dots, \un n\}$
 and
 $g\in S_n$ acts nontrivially only on additional random guests
 $\un{\ov 1}$,\dots, $\un{\ov n}$. These 'numbers' also are
 labels $\un 1$, \dots, $\un n$ on white triangles.
 We apply the construction of
 Subsect. \ref{ss:dense}
 and forget random guests.

\tt
Fakult\"at f\"ur Mathematik, Universit\"at Wien;

Institute for Information Transmission Problems, Moscow;

Dept. of Mechanics and Mathematics, Moscow State University.

e-mail:yurii.neretin@univie.ac.at

URL: https://mat.univie.ac.at/~neretin	
	
	\end{document}